\documentclass[12pt,twoside]{article}
\usepackage{times}
\usepackage{mathrsfs}
\usepackage{amsmath,amssymb}
\usepackage{color}
\usepackage{fancyhdr}
\allowdisplaybreaks

 \def \nn{\nonumber}

\newcommand{\be}{\begin{equation}}
\newcommand{\ee}{\end{equation}}
\newcommand{\bea}{\begin{eqnarray}}
\newcommand{\eea}{\end{eqnarray}}

\newcommand{\bR}{{\mathbb R}}
\newcommand{\bN}{{\mathbb N}}
\newcommand{\bZ}{{\mathbb Z}}

\def\p{\partial}
\def\la{\lambda}
\def\La{\Lambda}
\def\al{\alpha}
\def\t{\tilde}
\def\q{\quad}
\def\th{\theta}
\def\g{\gamma}

\def\dl{\delta}
\def\ve{\varepsilon}

\def\lt{\left}
\def\rt{\right}

\def\dl{\delta}
\def\Dl{\Delta}
\def\i{\infty}

\def \ls{\lesssim}
\def\p{\partial}
\def\f{\frac}
\def\na{\nabla}
\def\al{\alpha}

\def\t{\tilde}

\def\O{\Omega}

\def\q{\quad}
\def\qq{\qquad}
\def\s{\sqrt}

\def\B{\dot{\mathcal{B}}}
\def\E{\mathcal{E}}

\def\dB{\dot{B}}

\def\mX{\mathcal{X}}
\def\dD{\dot{\Delta}}
\def\ltr{\langle t\rangle^{\al}}

 \allowdisplaybreaks

\begin{document}
% \footskip=10pt
%\footnotesep=2pt
% \setcounter{page}{1}
\let\oldsection\section
\renewcommand\section{\setcounter{equation}{0}\oldsection}
\renewcommand\thesection{\arabic{section}}
\renewcommand\theequation{\thesection.\arabic{equation}}
\newtheorem{claim}{\noindent Claim}[section]
\newtheorem{theorem}{\noindent Theorem}[section]
\newtheorem{lemma}{\noindent Lemma}[section]
\newtheorem{proposition}{\noindent Proposition}[section]
\newtheorem{definition}{\noindent Definition}[section]
\newtheorem{remark}{\noindent Remark}[section]
\newtheorem{corollary}{\noindent Corollary}[section]
\newtheorem{example}{\noindent Example}[section]

\title{Global existence and optimal decay estimates of strong solutions to the compressible viscoelastic flows}

\author{Xinghong Pan ,\quad Jiang Xu
\vspace{0.5cm}
}

\date{}

\maketitle
% \vskip 0.2in

\centerline {\bf Abstract} \vskip 0.3 true cm
This paper is dedicated to the global existence and optimal decay estimates of strong solutions to the compressible viscoelastic flows in the whole space $\mathbb{R}^n$ with any $n\geq2$. We aim at extending those works by Qian \& Zhang and Hu \& Wang to the critical $L^p$ Besov space, which is not related to the usual energy space. With aid of intrinsic properties of viscoelastic fluids as in \cite{QZ1}, we consider a more complicated hyperbolic-parabolic system than usual Navier-Stokes equations. We define ``\emph{two effective velocities}'', which allows us to cancel out the coupling among the density, the velocity and the deformation tensor. Consequently, the global existence of strong solutions is constructed by using elementary energy approaches only. Besides, the optimal time-decay estimates of strong solutions will be shown in the general $L^p$ critical framework, which improves those decay results due to Hu \& Wu such that initial velocity could be \textit{large highly oscillating}.

\vskip 0.3 true cm

\vskip 0.3 true cm


\section{Introduction}

We consider the following equations of multi-dimensional compressible viscoelastic flows:

\be
\left\{
\begin{aligned}
&\p_t \rho+\textmd{div}(\rho u)=0, \\
&\p_t(\rho u)+ \textmd{div}(\rho u\otimes u)-\text{div}(2\mu D(u)+\la \text{div}u \text{Id})+\nabla P=\al\textmd{div}(\rho FF^T),\\
&\p_tF+u\cdot\nabla F=\nabla uF,
\end{aligned}
\right.\label{1.1}
\ee
where $\rho\in \bR_+$ is the density, $u\in\bR^n$ stands for the velocity and $F\in \bR^{n\times n}$ is the deformation gradient. $F^T$ means the transpose matrix of $F$. The pressure $P$ depends only upon the density and the function will be taking suitably smooth. Notations $\mathrm{div}$, $\otimes$ and $\nabla$
denote the divergence operator, Kronecker tensor product and gradient operator, respectively. $D(u)=\f{1}{2}(\na u+\na u^T)$ is the strain tensor. The density-dependent viscosity coefficients $\mu,\ \la$ are assumed to be smooth and to satisfy $\mu>0,\ \nu\triangleq \la+2\mu>0$. For simplicity,
the elastic energy $W(F)$ in system \eqref{1.1} has been taken to be the special
form of the Hookean linear elasticity:
\[
W(F)=\f{\al}{2}|F|^2, \q \al>0,
\]
which, however, does not reduce the essential difficulties for analysis. The methods and results of this
paper could be applied to more general cases.

In this paper, we focus on the Cauchy problem of system \eqref{1.1}, so the corresponding initial data are supplemented by
\be
(\rho,F;u)|_{t=0}=(\rho_0(x),F_0(x);u_0(x)), \q x\in\bR^n.  \label{id0}
\ee

It is well known that there are some fluids do not satisfy the classical Newtonian law. Also, there have been many attempts to capture different phenomena for non-Newtonian fluids, see for example \cite{ELZ1,Gme1,LLZ1,LW1} \textit{etc.}. System \eqref{1.1} is compressible viscoelastic flow of the Oldroyd type exhibiting the elastic behavior, which is one of the non-Newtonian fluids. We are interested in the well-posedness and stability of solutions to the Cauchy problem \eqref{1.1}-\eqref{id0}, at least under the perturbation of constant equilibrium state $(1,I,0)$.

Here let's first recall previous efforts related to viscoelastic flows. For the incompressible viscoelastic flows, there has been much
important progress on classical solutions. Lin-Liu-Zhang \cite{LLZ1}, Chen-Zhang \cite{CZ1}, Lei-Liu-Zhou\cite{LLZ2} and Lin-Zhang \cite{LZ1} established the local and global well-posedness with small data in Sobolev space $H^s$. Hu-Wu \cite{HW3} proved the long-time behavior and weak-strong uniqueness of solutions.
Chemin-Masmoudi \cite{CM1} proved the existence of a local solution and a global small solution in critical Besov spaces, where the Cauchy-Green strain tensor
is available in the evolution equation. Qian \cite{Qj1} proved the well-posedness of the incompressible viscoelastic system in critical spaces. Subsequently, Zhang-Fang \cite{ZF1} proved the global well-posedness in the critical $L^p$ Besov space. On the other hand, the global existence of weak solutions is still an open problem. Lions and Masmoudi \cite{LM1} considered a special case that the contribution of the strain rate is neglected, and proved the global existence of a weak solution with general initial data.

For compressible viscoelastic flows, Lei-Zhou \cite{LZ1} proved the global existence of classical solutions for the two-dimensional Oldroyd model via the incompressible limit. The local existence of strong solutions was obtained by Hu-Wang \cite{HW4}. Shortly, Hu-Wang \cite{HW1} and Qian-Zhang \cite{QZ1} independently proved the global existence in the critical $L^2$ Besov space with initial data near equilibrium. For convenience of reader, we would like to state their results as follows.

\begin{theorem}\label{thglobal1} (\cite{HW1,QZ1})
Assume that $P'(1)>0$. Then there exists two constant $\eta$ and $M$ such that if
\be
(\rho_0-1,F_0-I;u_0)\in\Big( \B^{n/2-1,n/2}_{2,2}\Big)^{1+n^2}\times \Big(\dot{B}^{n/2-1}_{2,1}\Big)^n. \nn
\ee
and
\be\label{1.3r}
\|(\rho_0-1,F_0-I)\|_{\B^{n/2-1,n/2}_{2,2}}+\|u_0\|_{\dot{B}^{n/2-1}_{2,1}}\leq \eta,\nn
\ee
then there exists a global unique solution $(\rho,F;u)$ of \eqref{1.1}$-$\eqref{id0} such that
\be
\|(\rho-1,F-I;u)\|_{\t{E}^{n/2}}\leq M\big(\|(\rho_0-1,F_0-I)\|_{\B^{n/2-1,n/2}_{2,2}}+\|u_0\|_{\dot{B}^{n/2-1}_{2,1}}\big),\nn
\ee
\end{theorem}
where
\bea
&&\t{E}^{n/2}\triangleq\Big(\t{\mathcal{C}_b}(\bR_+;\B^{n/2-1,n/2}_{2,2})\cap L^1(\bR_+;\B^{n/2+1,n/2}_{2,2})\Big)^{1+n^2}\times\nn\\
&&\qq\qq\Big(\t{\mathcal{C}_b}(\bR_+;\dot{B}^{n/2-1}_{2,1})\cap L^1(\bR_+;\dot{B}^{n/2+1}_{2,1})\Big)^{n},\nn
\eea
 with its norm, for $(a,O;v)\in \t{E}^{n/2}$,
\bea
\|(a,O;v)\|_{\t{E}^{n/2}}=\|(a,O)\|_{\t{L}^\i\B^{n/2-1,n/2}_{2,2}\cap L^1\B^{n/2+1,n/2}_{2,2}}+\|v\|_{\t{L}^\i\dB^{n/2-1}_{2,1}\cap L^1\dB^{n/2+1}_{2,1}}.\nn
\eea

Concerning those norm notations for the hybrid Besov space $\t{L}^q\B^{s,\sigma}_{2,p}(p\geq 2)$ and $\t{\mathcal{C}_b}(\B^{s,\sigma}_{2,p})$, the reader is referred to Section 2 below. Theorem \ref{thglobal1} stems from the scaling consideration. Note that system \eqref{1.1} is scaling invariant under the following transformation: for any constant $\kappa$,
\bea
&&(\rho_0(x),F_0(x);u_0(x))\rightarrow (\rho_0(\kappa x),F_0(\kappa x);\la u_0(\kappa x)),\nn\\
&&(\rho(t,x),F(t,x);u(t,x))\rightarrow (\rho(\kappa^2t, \kappa x),F(\kappa^2t,\kappa x);\kappa u(\kappa^2t,\kappa x)),\nn
\eea
up to changes of the Pressure $P$ into $\kappa^2P$ and the constant $\al$ into $\kappa^2\al$. This indicates the following definition of the critical Besov space.
\begin{definition}\label {defn1.1}
 A functional space is called a critical space if the associated norm is invariant under the transformation $(\rho(t,x),F(t,x);u(t,x))\rightarrow (\rho(\kappa^2t, \kappa x),F(\kappa^2t,\kappa x);\kappa u(\kappa^2t,\kappa x))$ (up to a constant independent of $\kappa$).
\end{definition}
Obviously, it is easy to see that  $\big(\dot{B}^{n/2}_{2,1}\big)^{1+n^2}\times\big(\dot{B}^{n/2-1}_{2,1}\big)^n$ is the critical space according to Definition \ref{defn1.1}. It should be emphasized that such basic idea is motivated by the seminal paper \cite{Dr1}, where the author first proved the global well-posedness for the compressible Navier-Stokes equations near equilibrium. Compared to \cite{Dr1}, there is an outstanding difficulty for the compressible viscoelastic system, that is, how to capture the damping effect of the deformation tensor among more complicated coupling between the velocity,
the density and the deformation tensor. Hu-Wang \cite{HW1} and Qian-Zhang \cite{QZ1} independently explored some intrinsic properties
of the viscoelastic system and established uniform estimate for more complicated linearized hyperbolic-parabolic systems, which eventually leads to Theorem \ref{thglobal1}.

The goal of this paper is twofold: firstly, we aim at extending the above statement (Theorem \ref{thglobal1}) to the critical $L^p$ Besov space, which allows highly large oscillating initial velocity. Secondly, we shall exhibit the long time behavior of the constructed solution.

Denote
\bea
\mathcal{E}^{n/p}&\triangleq& \Big\{(a,O;v)\big |(a,O;v)\in\Big(\t{\mathcal{C}_b}(\bR_+;\B^{n/2-1,n/p}_{2,p})\cap L^1(\bR_+;\B^{n/2+1,n/p}_{2,p})\Big)^{1+n^2}\nn\\
             &&\times\Big(\t{\mathcal{C}_b}(\bR_+;\B^{n/2-1,n/p-1}_{2,p})\cap L^1(\bR_+;\B^{n/2+1,n/p+1}_{2,p})\Big)^n\Big\},\nn
\eea
 with its norm
\be
\|(a,O;v)\|_{\E^{n/p}}=\|(a,O)\|_{\t{L}^\i\B^{n/2-1,n/p}_{2,p}\cap L^1\B^{n/2+1,n/p}_{2,p}}+\|v\|_{\t{L}^\i\B^{n/2-1,n/p-1}_{2,p}\cap L^1\B^{n/2+1,n/p+1}_{2,p}}.\nn
\ee
Now, we state the first result as follows.

\begin{theorem}\label{thglobal2}
Assume that $P'(1)>0$. Let $p$ satisfying $2\leq p\leq \min\big(4,2n/(n-2)\big)$ and, additionally, $p\neq 4$ if $n=2$.
If there exists two constant $\eta$ and $M$ such that if
\be
(\rho_0-1,F_0-I;u_0)\in\Big( \B^{n/2-1,n/p}_{2,p}\Big)^{1+n^2}\times \Big(\B^{n/2-1,n/p-1}_{2,p}\Big)^n. \nn
\ee
and
\be\label{1.3rx}
\|(\rho_0-1,F_0-I)\|_{\B^{n/2-1,n/p}_{2,p}}+\|u_0\|_{\B^{n/2-1,n/p-1}_{2,p}}\leq \eta,
\ee
then the Cauchy problem \eqref{1.1}-\eqref{id0} has a global unique solution $(\rho,F;u)$ such that $(\rho-1,F-I;u)
\in \mathcal{E}^{n/p}$ and
\be
\|(\rho-1,F-I;u)\|_{\mathcal{E}^{n/p}}\leq M\big(\|(\rho_0-1,F_0-I)\|_{\B^{n/2-1,n/p}_{2,p}}+\|u_0\|_{\B^{n/2-1,n/p-1}_{2,p}}\big).\nn
\ee
\end{theorem}

In comparison with those results in critical $L^p$ framework for compressible Navier-Stokes equations (see for example \cite{CD1} and \cite{CMZ1}), Theorem \ref{thglobal2} is not so surprising. Let's point out some new ingredients in the current proofs. To the best of our knowledge, there is a technical difficulty arising from a loss of one derivative for compressible N-S fluids, since there is no smoothing effect for the density in high frequency. To eliminate it, their proofs heavily rely on a paralinearized version combined with a Lagrangian change of variables, see \cite{CD1,CMZ1} for details. To the compressible viscoelastic system, the situation becomes more complicated. As shown by \cite{QZ1}, the damping effect of $F$ can be produced by some intrinsic conditions (see Proposition 3.1), however, similar to the density, there is not any smoothing effect at high frequencies. Here, in order to solve \eqref{1.1} globally, we follow an elementary energy approach in terms of \textit{effective velocity} rather than the elaborate Lagrangian change. The argument has been developed by
Haspot \cite{Hb1,Hb2} for compressible Navier-Stokes equations, which is based on the use of Hoff's viscous effective flux in \cite{Hd1}. Here, we introduce
the following ``\emph{two effective velocities}'',
$$w=\na(-\Dl)^{-1}(2a-\text{div} v),\qq \O^{ij}= e^{ij}+\f{1}{\mu_0}\La(-\Dl)^{-1} O^{ij}.$$
Indeed, the definition of $w$ is almost the same as that in \cite{Hb1,Hb2}. A slight difference lies on the coefficient of $a$, which comes from contribution of the deformation gradient $F$. Another effective velocity with respect to $\O^{ij}$ is new, which allows to cancel the coupling between $e^{ij}$ and $O^{ij}$ at high frequencies (see Sections 4 and 5 for more details). In physical dimensions $n=2,3,$  the value of $p$ enable us to consider the case  $p>n$ for which the
velocity regularity exponent $n/p-1$  becomes negative. Consequently, Theorem \ref{thglobal2} applies  to \emph{large} highly oscillating initial velocities
(see \cite{CD1,CMZ1} for more explanation).

An interesting question follows after gaining Theorem \ref{thglobal2}. One may wonder how the global strong solutions constructed above look like for large
time. Although providing an accurate long-time asymptotic description is still out of reach, a number of results concerning the time decay rates of global solutions, sometimes referred to as $L^q-L^r$ decay rates are available. For example, Hu-Wu \cite{HW2} proved the global existence of strong solutions to \eqref{1.1} as initial data are the small perturbation $(1,I;0)$ in $H^2(\bR^3)$. Furthermore, with the extra $L^1(\bR^3)$ assumption, it was shown that
those solutions converged to equilibrium state at the following way
\be
\|(\rho-1,F-I;u)\|_{L^p}\leq C\langle t\rangle^{-\f{3}{2}(1-\f{1}{p})}. \label{1.3}
\ee
The decay rate in \eqref{1.3} turns out the same one for the heat kernel, which is sometime referred as the optimal time-decay rate. Next, we state a decay result for those solutions constructed in  Theorem \ref{thglobal2}. Precisely, one has
\begin{theorem}\label{thdecay}
Let $n\geq2$ and $p$ satisfies $2\leq p\leq \min\big(4,2n/(n-2)\big)$ and $p\neq 4$ if $n=2$. Let $(\rho_0,u_0,F_0)$ fulfill the assumptions of Theorem\ref{thglobal2} and $(\rho,u,F)$ be the global solution of System \eqref{1.1}. Then there exists a constant $\sigma=\sigma(p,n,\la,\mu,\al,P)$ such that if additionally
\be
\mathcal{G}_{p,0}\triangleq\|(\rho_0-1,F_0-I;u_0)\|^\ell_{\dot{B}^{-s_0}_{2,\i}}\leq \sigma \q with \q s_0\triangleq n(2/p-1/2), \label{id3}
\ee
then we have for $t\geq0$,
\be \label{de1}
\mathcal{G}_p(t)\ls\Big(\mathcal{G}_{p,0}+\|(\nabla \rho_0,\na F_0;u_0)\|^h_{\dot{B}^{n/p-1}_{p,1}}\Big),
\ee
where $\mathcal{G}_p(t)$ is defined by
\bea \label{dfD}
&&\mathcal{G}_p(t)\triangleq\sup\limits_{s\in[\ve-s_0,\f{n}{2}+1]}\|\langle\tau\rangle^{\f{s_0+s}{2}}(\rho-1,F-I;u)\|^\ell_{\t{L}^\i_t\dot{B}^s_{2,1}}+\|\langle\tau\rangle^{\al}(\nabla a,\nabla F;u)\|^h_{\t{L}^\i_t\dot{B}^{\f{n}{p}-1}_{p,1}}\nn\\
                 &&\qq\qq+\|\tau \nabla u\|^h_{\t{L}^\i_t\dot{B}^{\f{n}{p}}_{p,1}},
\eea
with $\al:=n/p+1/2-\ve$ and $\ve>0$ is a sufficiently small constant.
\end{theorem}

Here $\|f\|^\ell_{\bullet}$ and $\|f\|^h_{\bullet}$ represent the low and high frequency part of some norm $\|f\|_{\bullet}$ to a tempered distribution $f$ whose exact definition will be given in Section 2.

Some comments are in order.
\begin{enumerate}

\item Due to the Sobolev imbedding properties $L^1\hookrightarrow \dot{B}^{0}_{1,\infty}\hookrightarrow\dot{B}^{-n/2}_{2,\infty},\ \ \dot{H}^{-n/2}\hookrightarrow \dot{B}^{-n/2}_{2,\infty}$, our low-frequency assumption is less restrictive. Actually,
the assumption is relevant in other contexts like the Boltzmann equation (see \cite{SS1}), or hyperbolic systems with dissipation (see \cite{XK2}).
\item The decay result remains true in the case of \emph{large} highly oscillating initial velocities, since the case $p>n$ occurs in physical dimensions $n=2,3,$ which was not shown by previous efforts (see for example \cite{HW2}).
\item Likewise,  ``\emph{two effective velocities}'' play a key role in establishing the nonlinear time-weighted inequality (\ref{dfD}). Furthermore, the
optimal decay estimates of  $L^q$-$L^r$ type can be derived from the definition of $\mathcal{G}_p(t)$ by using interpolation tricks. The interested reader is referred to \cite{DX1} for similar details.
\end{enumerate}

The rest of this paper is arranged as follows: In Section 2, we first review the Littlewood-Paley theory and give definitions and estimates for the hybrid-Besov space. In Section 3, we reformulate our system into a hyperbolic-parabolic system coupled by the density, the velocity and the deformation gradient. Section 4 is devoted to presenting the proof of Theorem\ref{thglobal2}. In Section 5, we prove the decay estimate in Theorem \ref{thdecay}.
Some analysis properties in the hybrid Besov space are also given in the Appendix.

\section{Littlewood-Paley Theory and the Hybrid Besov Space }
\q\ Throughout the paper, we denote by $C$ a generic constant which may be different from line to line. The notation $A\ls B$ means $A\leq CB$ and $A\thickapprox B$ indicates $A\leq CB$ and $B\leq CA$.
\subsection{Littlewood-Paley decomposition}

\q\ Let's begin with the Littlewood-Paley decomposition. There exists two radial smooth functions $\varphi(x),\chi(x)$ supported in the annulus $\mathcal{C}=\{\xi\in\bR^n:3/4\leq |\xi|\leq 8/3\}$ and the ball $B=\{\xi\in\bR^n:|\xi|\leq 4/3\}$, respectively such that
\be
\chi(\xi)+\sum\limits_{j\geq 0}\varphi(2^{-j}\xi)=1\q \forall \xi\in\bR^n.\nn
\ee
\be
\sum\limits_{j\in \bZ}\varphi(2^{-j}\xi)=1\q \forall \xi\in\bR^n\diagdown\{0\}.\nn
\ee

The homogeneous dyadic blocks $\dD_j$ and the homegeneous low-frequency cut-off operators $\dot{S}_j$ are defined for all $j\in \bZ$ by
\be
\dD_j u=\varphi(2^{-j}D)f,\q \dot{S}_jf=\sum\limits_{k\leq j-1}\dD_k f=\chi(2^{-j}D)f.\nn
\ee

The following Bernstein inequality will be repeatedly used throughout the paper.

\begin{lemma}[\cite{BCD1}]\label{lem2.1}
Let $\mathcal{C}$ be an annulus and $B$ a ball. A constant $C$ exists such that for any nonnegative integer $k$, any couple $(p,q)$ in $[1,\i]^2$ with $q\geq p\geq 1$, and any function $u$ of $L^p$, we have
\be
\text{Supp} \hat{u}\subset\la B\Rightarrow \sup\limits_{|\al|=k}\|\p^\al u\|_{L^q}\leq C^{k+1}\la^{k+n(\f{1}{p}-\f{1}{q})}\|u\|_{L^p}, \nn
\ee
\be
\text{Supp} \hat{u}\subset\la \mathcal{C}\Rightarrow C^{-k-1}\la^k\|u\|_{L^p}\leq\sup\limits_{|\al|=k}\|\p^\al u\|_{L^p}\leq C^{k+1}\la^{k}\|u\|_{L^p}. \nn
\ee
\end{lemma}
\subsection{The hybrid Besov space}

\q\ We denote by\ $\mathcal{Z}'(\bR^n)$ the dual space of
\be
\mathcal{Z}(\bR^n)\triangleq \{f\in \mathcal{S}(\bR^n):\p^\al \hat{f}(0)=0,\forall \al\in(\bN\cup 0)^n\}.\nn
\ee
Firstly, we give the definition of the the homogeneous Besov space.
\begin{definition}
Let s be a real number and (p,r) be in $[1,\i]^2$. The homogeneous Besov space $\dot{B}^s_{p,r}$ consists of those distributions $u\in \mathcal{Z}'(\bR^n)$ such that
\be
\|u\|_{\dot{B}^s_{p,r}}\triangleq \Big(\sum\limits_{j\in\bZ}2^{jsr}\|\dD_j u\|^r_{L^p}\Big)^{\f{1}{r}}<\i.\nn
\ee
\end{definition}

Secondly, we introduce the hybrid Besov space that will be used in this paper.
\begin{definition}
Let $s,\sigma\in \bR$, $1\leq p\leq +\i$. The hybrid Besov space $\B^{s,\sigma}_{2,p}$ is defined by
\be
\B^{s,\sigma}_{2,p}\triangleq \{f\in\mathcal{Z}'(\bR^n):\|f\|_{\B^{s,\sigma}_{2,p}}<\i\}, \nn
\ee
with
\be
\|f\|_{\B^{s,\sigma}_{2,p}}\triangleq \sum\limits_{2^k\leq R_0}2^{ks}\|\dD_k f \|_{L^2}+\sum\limits_{2^k>R_0}2^{k\sigma}\|\dD_k f \|_{L^p},\nn
\ee
\end{definition}
where $R_0$ is a fixed and sufficiently large constant which may depending on $\la(1),\mu(1),p$ and $n$.

Since we are concerned with time-dependent functions valued in Besov spaces, the space-time mixed norm is usually given by
\be
\|u\|_{L^q_T\B^{s,\sigma}_{2,p}}:=\big\|\|u(t,\cdot)\|_{\B^{s,\sigma}_{2,p}}\big\|_{L^q(0,T)}.\nn
\ee
Here, we introduce another space-time mixed Besov norm, which is referred to Chemin-Lerner's spaces. The definition is as follows.
\be
\|u\|_{\t{L}^q_T\dot{B}^{s,\sigma}_{2,p}}\triangleq\sum\limits_{2^k\leq R_0}2^{ks}\|\dD_k u\|_{L^q(0,T)L^2}+\sum\limits_{2^k> R_0}2^{k\sigma}\|\dD_k u\|_{L^q(0,T)L^2}.\nn
\ee
The index $T$ will be omitted if $T=+\i$ and we shall denote by $\t{\mathcal{C}}_b(\dB^{s,\sigma}_{2,p})$ the subset of functions $\t{L}^\i(\dB^{s,\sigma}_{2,p})$ which are continuous from $\bR_+$ to $\dB^{s,\sigma}_{2,p}$. It is easy to check that $\t{L}^1_T\B^{s,\sigma}_{2,p}=L^1_T\B^{s,\sigma}_{2,p}$ and $\t{L}^q_T\B^{s,\sigma}_{2,p}\subseteq L^q_T\B^{s,\sigma}_{2,p}$ for $q>1$.

Also, for a tempered distribution $f$, we denote
\be
f^\ell\triangleq \sum\limits_{2^k\leq R_0}\dD_k f,\q f^h\triangleq f-f^\ell,\nn
\ee
and
\be
\|f\|^\ell_{\dB^{s}_{p,1}}=\sum\limits_{2^k\leq R_0}2^{ks}\|\dD_k f \|_{L^p},\q \|f\|^h_{\dB^{s}_{p,1}}=\sum\limits_{2^k>R_0}2^{k\sigma}\|\dD_k f \|_{L^p},\nn
\ee

\be
\|f\|^\ell_{\t{L}^q_T\dB^{s}_{p,1}}=\sum\limits_{2^k\leq R_0}2^{ks}\|\dD_k f \|_{L^q(0,T;L^p)},\q \|f\|^h_{\t{L}^q_T\dB^{s}_{p,1}}=\sum\limits_{2^k>R_0}2^{k\sigma}\|\dD_k f \|_{L^q(0,T;L^p)},\nn
\ee
\be
\|f\|^\ell_{\dB^{s}_{2,\i}}=\sup_{2^k\leq R_0}\limits2^{ks}\|\dD_k f \|_{L^2},\nn
\ee
for $s\in \bR$.

Next, we collect nonlinear estimates in Besov spaces.
\begin{lemma}\label{l2.2}
For the Besov space, we have the following properties:

1,$\B^{s_2,\sigma}_{2,p}\subseteq \B^{s_1,\sigma}_{2,p}$ for $s_1\geq s_2$ and $\B^{s,\sigma_2}_{2,p}\subseteq \B^{s,\sigma_1}_{2,p}$ for $\sigma_1\leq \sigma_2$.

2,Interpolation: For $s_1,s_2,\sigma_1,\sigma_2\in \bR$ and $\theta\in [0,1]$, we have
\be
\|f\|_{\B^{\theta s_1+(1-\theta)s_2,\theta\sigma_1+(1-\theta)\sigma_2}_{2,p}}\leq \|f\|^\theta_{\B^{s_1,\sigma_1}_{2,p}}\|f\|^{(1-\theta)}_{\B^{s_2,\sigma_2}_{2,p}}.\nn
\ee

3,Embedding: $L^\i\hookrightarrow \B^{n/2,n/p}_{2,p}$;\\

\qq\qq\q\q $\dot{B}^s_{2,1}\hookrightarrow \B^{s,s-n/2+n/p}_{2,p}\hookrightarrow \dot{B}^{s-n/2+n/p}_{p,1}$ for $p\geq 2$.
\end{lemma}
\begin{lemma}[\cite{Dr2}\label{ll2.3}]
Let $1\leq p,q,q_1,q_2\leq \i$ with $\f{1}{q_1}+\f{1}{q_2}=\f{1}{q}$. Then we have the following:

(1),If $s_1,s_2\leq n/p$ and $s_1+s_2>n \max(0,2/p-1)$, then
\be
\|fg\|_{\t{L}^q_T(\dot{B}^{s_1+s_2-n/p}_{p,1})}\leq C\|f\|_{\t{L}^{q_1}_T(\dot{B}^{s_1}_{p,1})}\|g\|_{\t{L}^{q_2}_T(\dot{B}^{s_2}_{p,1})}.\nn
\ee

(2),If $s_1\leq n/p,\ s_2<n/p$ and $s_1+s_2>n \max(0,2/p-1)$, then
\be
\|fg\|_{\t{L}^q_T\dot{B}^{s_1+s_2-n/p}_{p,\i}}\leq C\|f\|_{\t{L}^{q_1}_T(\dot{B}^{s_1}_{p,1})}\|g\|_{\t{L}^{q_2}_T(\dot{B}^{s_2}_{p,\i})}.\nn
\ee
\end{lemma}
\begin{remark}
Lemma\ref{ll2.3} still remain true in the usual homogenous Besov spaces. For example the estimate in Lemma\ref{ll2.3}(1) becomes
\be
\|fg\|_{\dot{B}^{s_1+s_2-n/p}_{p,1}}\leq C\|f\|_{\dot{B}^{s_1}_{p,1}}\|g\|_{\dot{B}^{s_2}_{p,1}}.\nn
\ee
\end{remark}

\begin{lemma}
[\cite{DX1}\label{l2.4}] Let $\sigma>0$ and $1\leq p,r\leq \i$. Then $\dB^\sigma_{p,r}\cap L^\i$ is an algebra and
\be
\|fg\|_{\dB^{\sigma}_{p,r}}\ls \|f\|_{L^\i}\|g\|_{\dB^{\sigma}_{p,r}}+\|g\|_{L^\i}\|f\|_{\dB^{\sigma}_{p,r}}.\nn
\ee
Let $\sigma_1,\sigma_2,p_1,p_2$ satisfy
\be
\sigma_1+\sigma_2>0,\ \sigma_1\leq n/p_1,\ \sigma_2\leq n/p_2,\ \sigma_1\geq\sigma_2,\ \f{1}{p_1}+\f{1}{p_2}\leq 1.\nn
\ee
Then we have
\be\label{ee2.1}
\|fg\|_{\dB^{\sigma_2}_{q,1}}\ls \|f\|_{\dB^{\sigma_1}_{p_1,1}}\|f\|_{\dB^{\sigma_2}_{p_2,1}} \q \text{with}\ \f{1}{q}=\f{1}{p_1}+\f{1}{p_2}-\f{\sigma_1}{n}.
\ee
Finally for exponents $\sigma>0,\ 1\leq p_1,p_2,q\leq \i$ satisfying
\be
\f{n}{p_1}+\f{n}{p_2}-n\leq \sigma\leq\min(\f{n}{p_1},\f{n}{p_2})\q \text{and}\ \f{1}{q}=\f{1}{p_1}+\f{1}{p_2}-\f{\sigma}{n},\nn
\ee
we have
\be\label{ee2.2}
\|fg\|_{\dB^{-\sigma}_{q,\i}}\ls\|f\|_{\dB^{\sigma}_{p_1,1}}\|g\|_{\dB^{-\sigma}_{p_2,\i}}.
\ee
\end{lemma}

\begin{lemma}
[\cite{DX1}\label{l2.5}]There exists a universal interger $N_0$ such that for any $2\leq p\leq 4$, and $\sigma>0$, we have
\be\label{e2.1}
\|fg^h\|^\ell_{\dB^{-s_0}_{2,\i}}\ls\big(\|f\|_{\dB^\sigma_{p,1}}+\|\dot{S}_{k_0+N_0}f\|_{L^{p^\ast}}\big)\|g^h\|_{\dB^{-\sigma}_{p,\i}},
\ee
\be\label{e2.2}
\|f^hg\|^\ell_{\dB^{-s_0}_{2,\i}}\ls\big(\|f^h\|_{\dB^\sigma_{p,1}}+\|\dot{S}_{k_0+N_0}f^h\|_{L^{p^\ast}}\big)\|g\|_{\dB^{-\sigma}_{p,\i}},
\ee
with $s_0= n\big(\f{2}{p}-\f{1}{2}\big)$ and $\f{1}{p^\ast}=\f{1}{2}-\f{1}{p}$.
\end{lemma}
\begin{lemma}
[\cite{DX1}\label{ll2.6}] Let $1\leq p,p_1\leq \i$ and
\be
-\min\big(\f{n}{p_1},\f{n}{p'}\big)<\sigma\leq 1+\min\big(\f{n}{p},\f{n}{p_1}\big).\nn
\ee
There exists a constant $C>0$, depending only on $\sigma$ such that for all $j\in\bZ$, we have
\be\label{e2.3}
\|[v\cdot\na,\na \dD_j]z\|_{L^p}\leq Cc_j2^{-j(\sigma-1)}\|\na v\|_{\dB^{n/p_1}_{p_1,1}}\|\na z\|_{\dB^{\sigma-1}_{p,1}},
\ee
where $(c_j)_{j\in\bZ}$ denotes a sequence such that $\|(c_j)\|_{\ell^1}\leq 1$.
\end{lemma}

\section{Reformulation of System \eqref{1.1}}

\q\ Here, we present intrinsic properties of compressible viscoelastic flows, which have been explored in \cite{QZ1}.
\begin{proposition}\label{ipp}
The density $\rho$ and the deformation gradient $F$ of \eqref{1.1} satisfy the following relations:
\be
\nabla\cdot(\rho F^T)=0 \q\text{and}\q F^{lk}\p_lF^{ij}-F^{lj}\p_lF^{ik}=0, \label{ip1}
\ee
if the initial data $(\rho_0,F_0)$ satisfies
\be
\nabla\cdot(\rho_0 F_0^T)=0 \q\text{and}\q F_0^{lk}\p_lF_0^{ij}-F_0^{lj}\p_lF_0^{ik}=0. \label{ip2}
\ee
\end{proposition}

By Proposition \ref{ipp}, the $i$-th component of the vector $\text{div}(\rho FF^T)$ can be written as
\bea\label{e1}
\p_j(\rho F^{ik}F^{jk})&=&\rho F^{jk}\p_jF^{ik}+F^{ik}\p_j(\rho F^{jk})\nn\\
&=&\rho F^{jk}\p_jF^{ik},
\eea
where we used the first equality in \eqref{ip1}.

Denote $\chi_0=(P'(1))^{-1/2}$ and define
\be
a(t,x)=\rho(\chi^2_0t,\chi_0x)-1,\q v(t,x)=\chi_0u(\chi^2_0t,\chi_0x),\q O(t,x)=F(\chi^2_0t,\chi_0x)-I.\nn
\ee
By using \eqref{e1}, we get
\be\label{r1}
\left\{
\begin{aligned}
&\p_t a+v\cdot\nabla a+\nabla\cdot v=-a\nabla\cdot v, \\
&\p_tv+v\cdot\nabla v-\mathcal{A}v+\na a-\beta\na\cdot O=\beta O^{jk}\p_jO^{\bullet k}-I(a)\mathcal{A}v-K(a)\na a\\
&\qq\qq\qq\qq\qq\qq +\f{1}{1+a}\text{div}\big(2\t{\mu}(a)D(v)+\t{\la}(a)\text{div} v\text{Id}\big),\\
&\p_tO+v\cdot\nabla O-\nabla v=\nabla vO,
\end{aligned}
\right.
\ee
where
\be
I(a)\triangleq\f{a}{1+a},\ K(a)\triangleq \f{P'(1+a)}{(1+a)P'(1)}-1,\ \mathcal{A}=\mu(1)\Dl+(\la(1)+\mu(1))\nabla\text{div},
\nn
\ee
and
\be\q \beta=\f{\al}{P'(1)}, \ \t{\mu}(a)=\mu(1+a)-\mu(1),\  \t{\la}(a)=\la(1+a)-\la(1).\nn
\ee
$\ O^{jk}\p_jO^{\bullet k}$ is a vector function whose components are $(O^{jk}\p_jO^{i k})^n_{i=1}$. For simplicity, we set $\la(1)=\la_0, \mu(1)=\mu_0 $. Here and below, we normalize $\beta=1$ and $\nu(1):=\la(1)+2\mu(1)=1$ without loss of generality.

For $s\in \bR$, we denote
\be
\Lambda^sf\triangleq\mathcal{F}^{-1}(|\xi|^s\mathcal{F}(f)),\nn
\ee
and introduce two variables as in  \cite{QZ1}:
\be\label{3.5ee}
d=\Lambda^{-1}\text{div}v, \q e^{ij}=\La^{-1}\p_j v^i.
\ee
Using the second equality in \eqref{ip1}, we have
\be\label{3.6e}
\La^{-1} (\p_j\p_k O^{ik})=-\La O^{ij}-\La^{-1}\p_k(O^{lj}\p_lO^{ik}-O^{lk}\p_lO^{ij}).
\ee
Hence, with aid of \eqref{3.6e}, the system \eqref{r1} can be reformulated as follows
\be\label{r2}
\left\{
\begin{aligned}
&\p_t a+\Lambda d=G_1, \\
&\p_t e^{ij}-\mu_0\Dl e^{ij}-(\la_0+\mu_0)\p_i\p_j d+\La^{-1}\p_i\p_j a+\La O^{ij}=G^{ij}_4\\
&\p_tO^{ij}-\La e^{ij}=G^{ij}_3,\\
&d=-\La^{-2}\p_i\p_j e^{ij}, v^i=-\La^{-1}\p_je^{ij},
\end{aligned}
\right.
\ee
where $G_1=-a\nabla\cdot v-v\cdot\na a$, $G^{ij}_3=\p_k v^iO^{kj}-v\cdot\na O^{ij}$ and
\bea
&&G^{ij}_4=-\Lambda^{-1}\p_j\Big(v\cdot\nabla v^i-O^{lk}\p_lO^{i k}+I(a)(\mathcal{A}v)^i+K(a)\p_i a\Big)\nn\\
&&\qq\q-\La^{-1}\p_k(O^{lj}\p_lO^{ik}-O^{lk}\p_lO^{ij})\nn\\
&&\qq\q+\La^{-1}\p_j\Big(\f{1}{1+a}\text{div}\big(2\t{\mu}(a)D(v)+\t{\la}(a)\text{div} v\text{Id}\big)\Big)^i.\nn
\eea
Additionally, we need the auxiliary equation in subsequent estimates
\be\label{3.11e}
\p_i O^{ij}=-\p_ja-G^j_0,\qq G^j_0=\p_i(aO^{ij}),
\ee
 which is deduced from the first equality in \eqref{ip1}.
\section{Proof of Theorem \ref{thglobal2}}
Inspired by \cite{CMZ1}, we may extend those results in \cite{QZ1} to the $L^p$ critical framework.
First of all, it is convenient to give the following interpolation inequalities
\be\label{4.1er}
\begin{aligned}
&\|f\|_{\t{L}^2_T\B^{n/2,n/p}_{2,p}}\ls \|f\|^{1/2}_{\t{L}^{\i}_T\B^{n/2-1,n/p}_{2,p}}\|f\|^{1/2}_{{L}^{1}_T\dB^{n/2+1,n/p}_{2,p}},\\
&\|f\|_{\t{L}^2_T\B^{n/2,n/p}_{2,p}}\ls \|f\|^{1/2}_{\t{L}^{\i}_T\B^{n/2-1,n/p-1}_{2,p}}\|f\|^{1/2}_{{L}^{1}_T\dB^{n/2+1,n/p+1}_{2,p}}.
\end{aligned}
\ee
The proof Theorem \ref{thglobal2} is divided into several parts. The first one is to establish two a priori estimates.
\subsection{Two a priori estimates}

\q\ Let $T>0$. We denote the following functional space $\mathcal{E}^{n/p}_T$ by
\bea
\mathcal{E}^{n/p}_T&\triangleq& \big\{(a,O;v)\in\big(\t{L}^\i(0,T;\B^{n/2-1,n/p}_{2,p})\cap L^1(0,T;\B^{n/2+1,n/p}_{2,p})\big)^{1+n^2}\nn\\
             &&\times\big(\t{L}^\i(0,T;\B^{n/2-1,n/p-1}_{2,p})\cap L^1(0,T;\B^{n/2+1,n/p+1}_{2,p})\big)^{n}\big\}\nn
\eea
with the norm
\bea
\|(a,O;v)\|_{\mathcal{E}^{n/p}_T}&\triangleq& \|(a,O)\|_{\t{L}^\i_TB^{n/2-1,n/p}_{2,p}\cap L^1_T\B^{n/2+1,n/p}_{2,p}}+\|v\|_{\t{L}^\i_T\B^{n/2-1,n/p-1}_{2,p}\cap L^1_T\B^{n/2+1,n/p+1}_{2,p}}.\nn
\eea
\begin{proposition}\label{pro4.1}
Let $2\leq p\leq \min(4,\f{2n}{n-2})$ and $p<2n$. Assume that $(a,O;v)$ is a strong solution of system \eqref{r1} on $[0,T]$ with
\be
\|a\|_{L^\i([0,T]\times\bR^n)}\leq \f{1}{2}.  \nn
\ee
Then we have
\bea\label{4.1}
&&\|(a,O;v)\|_{\E^{n/p}_T}\leq C\Big\{\|(a_0,O_0;v_0)\|_{\E^{n/p}_0}\nn\\
&&\qq\qq\qq\qq+\|(a,O;v)\|^{2}_{\E^{n/p}_T}\big(1+\|(a,O;v)\|_{\E^{n/p}_T}\big)^{n+3}\Big\},
\eea
where $\|(a_0,O_0;v_0)\|_{\E^{n/p}_0}\triangleq \|(a_0,O_0)\|_{\B^{n/2-1,n/p}_{2,p}}+\|v_0\|_{\B^{n/2-1,n/p-1}_{2,p}}$.
\end{proposition}

We introduce another functional space $E^{n/2}_T$ defined by
\bea
&&E^{n/2}_T\triangleq \Big\{(a,O;v)\in\big(\t{L}^\i(0,T;\B^{n/2-1,n/2}_{2,2})\cap L^1(0,T;\B^{n/2+1,n/2}_{2,2})\big)^{1+n^2}\nn\\
&&\qq\qq\qq\qq\qq\qq \times\big(\t{L}^\i(0,T;\dot{B}^{n/2-1}_{2,1})\cap L^1(0,T;\dot{B}^{n/2+1}_{2,1})\big)^{n}\Big\}\nn
\eea
with the norm
\be
\|(a,O;v)\|_{E^{n/2}_T}\triangleq \|(a,O)\|_{\t{L}^\i_T\B^{n/2-1,n/2}_{2,2}\cap L^1_T\B^{n/2+1,n/2}_{2,2}}+\|v\|_{\t{L}^\i_T\dot{B}^{n/2-1}_{2,1}\cap L^1_T\dot{B}^{n/2+1}_{2,1}}.\nn
\ee

\begin{proposition}\label{pro4.2}
Under the assumption of Proposition \ref{pro4.1}, we have
\bea\label{4.8r}
&&\|(a,O;v)\|_{E^{n/2}_T}\leq C\Big\{\|(a_0,O_0;v_0)\|_{E^{n/2}_0}\nn\\
&&\qq\qq\qq\qq+\|(a,O;v)\|_{E^{n/2}_T}\|(a,O;v)\|_{\E^{n/p}_T}\big(1+\|(a,O;v)\|_{\E^{n/p}_T}\big)^{n+3}\Big\},\nn\\
\eea
where $\|(a_0,O_0;v_0)\|_{E^{n/2}_0}\triangleq \|(a_0,O_0)\|_{\B^{n/2-1,n/2}_{2,2}}+\|v_0\|_{\dot{B}^{n/2-1}_{2,1}}$.
\end{proposition}

The proof of Propositions \ref{pro4.1}-\ref{pro4.2} lie in the pure energy methods in terms of low-frequency and high-frequency decompositions.\\

\noindent\textbf{Step1: Low-frequency estimates ($2^k\leq R_0$).}\\
Denote $a_k=\dD_k a, O_k=\dD_k O$ and $d_k=\dD_k d, e_k=\dD_k e$ for simplicity. By applying $\dD_k$ to \eqref{r2}, we have
\be\label{e4.3}
\left\{
\begin{aligned}
&\p_t a_k+\Lambda d_k=\dD_k G_1, \\
&\p_t e^{ij}_k-\mu_0\Dl e^{ij}_k-(\la_0+\mu_0)\p_i\p_j d_k+\La^{-1}\p_i\p_j a_k+\La O^{ij}_k=\dD_k G^{ij}_4\\
&\p_tO^{ij}_k-\La e^{ij}_k=\dD_kG^{ij}_3,\\
&d_k=-\La^{-2}\p_i\p_j e^{ij}_k.
\end{aligned}
\right.
\ee
Taking $L^2$ inner product of $\eqref{e4.3}_2$ with $e^{ij}_k$, and then summing up the resulting equation with respect to indices $i,j$, we arrive at
\bea\label{e4.4}\f{1}{2}\|e_k\|^2_{L^2}+\mu_0\|\La e_k\|^2_{L^2}+(\la_0+\mu_0)\|\La d_k\|^2_{L^2}-(a_k|\La d_k)+(\La O_k|e_k)=(\dD_k G_4|e_k),
\eea
where we have used the fact $d_k=-\La^{-2}\p_i\p_j e^{ij}_k$.

 Taking $L^2$ inner product of $\eqref{e4.3}_1$ and $\eqref{e4.3}_3$ with $a_k$ and $O_k$, respectively, and then adding  the resulting equations to \eqref{e4.4} together, we obtain
\bea\label{e4.5}
&&\f{1}{2}\Big(\|a_k\|^2_{L^2}+\|O_k\|^2_{L^2}+\|e_k\|^2_{L^2}\Big)+\mu_0\|\La e_k\|^2_{L^2}+(\la_0+\mu_0)\|\La d_k\|^2_{L^2}\nn\\
&=&(\dD_k G_1|a_k)+(\dD_k G_4|e_k)+(\dD_k G_3|O_k).
\eea

To capture the dissipation arising from $(a,O)$, we next apply the operator $\La$ to $\eqref{e4.3}_1$ and take the $L^2$ inner product of the resulting equation with $-d_k$. Also, we take the $L^2$ inner product of $\eqref{e4.3}_2$ with $\La^{-1}\p_i\p_j a_k$. Therefore, we add those resulting equations and get
\bea\label{e4.6}
&&-\f{d}{dt}(\La a_k|d_k)+\|\La a_k\|^2_{L^2}-\|\La d_k\|^2_{L^2}-(\La^2 d_k|\La a_k)+(O^{ij}_k|\p_i\p_j a_k)\nn\\
&=&-(\La \dD_k G_1|d_k)+(\dD_kG^{ij}_4|\La^{-1}\p_i\p_ja_k).
\eea
On the other hand,  we apply $\La$ to $\eqref{e4.3}_3$ and then take the $L^2$ inner product of the resulting equation with $e^{ij}_k$. We also take the $L^2$ inner product of $\eqref{e4.3}_2$ with $\La O^{ij}_k$. By summing up those resulting equations, we obtain
\bea\label{e4.7}
&&\f{d}{dt}(\La O_k|e_k)+\|\La O_k\|^2_{L^2}-\|\La e_k\|^2_{L^2}\nn\\
&&-(\la_0+\mu_0)(\La O^{ij}_k|\p_i\p_j d_k)+\mu_0(\La^2e_k|\La O_k)+(\p_i\p_ja_k|O^{ij}_k)\nn\\
&=&(\La \dD_k G_3|e_k)+(\dD_kG_4|\La O_k).
\eea

Now, we multiply a small constant $\nu_1>0$ (to be determined) to \eqref{e4.6} and \eqref{e4.7}, respectively, and then add
the resulting equations with \eqref{e4.5} together. Consequently, we are led to the following inequality
\bea\label{e4.8}
&&\f{1}{2}\Big(\|a_k\|^2_{L^2}+\|O_k\|^2_{L^2}+\|e_k\|^2_{L^2}+2\nu_1(\La O_k|e_k)-2\nu_1(\La a_k|d_k)\Big)\nn\\
&&+(\mu_0-\nu_1)\|\La e_k\|^2_{L^2}+(\la_0+\mu_0-\nu_1)\|\La d_k\|^2_{L^2}+\nu_1(\|\La a_k\|^2_{L^2}+\|\La O_k\|^2_{L^2})\nn\\
&&+\nu_1\mu_0(\La^2e_k|\La O_k)-\nu_1(\la_0+\mu_0)(\La O^{ij}_k|\p_i\p_j d_k)-\nu_1(\La^2 d_k|\La a_k)+2\nu_1(\p_i\p_ja_k|O^{ij}_k)\nn\\
&=&(\dD_k G_1|a_k)+(\dD_k G_4|e_k)+(\dD_k G_3|O_k)-\nu_1(\La \dD_k G_1|d_k)\nn\\
&&+\nu_1(\dD_kG^{ij}_4|\La^{-1}\p_i\p_ja_k)+\nu_1(\La \dD_k G_3|e_k)+\nu_1(\dD_kG_4|\La O_k).
\eea
It follows from \eqref{3.11e} that
\bea\label{e4.9}
(\p_i\p_ja_k|O^{ij}_k)&=&(a_k|\p_i\p_jO^{ij}_k)\nn\\
  &=&\big((-\Dl a_k-\p_j\dD_k G^j_0)|a_k\big)\nn\\
  &=&\|\La a_k\|^2_{L^2}-( a_k|\p_j \dD_k G^j_0).
\eea
Inserting \eqref{e4.9} into \eqref{e4.8}, we can get
\bea\label{e4.10}
&&\f{d}{dt}f^2_{\ell,k}+\t{f}^2_{\ell,k}\nn\\
&=&(\dD_k G_1|a_k)+(\dD_k G_4|e_k)+(\dD_k G_3|O_k)\nn\\
&&-\nu_1(\La \dD_k G_1|d_k)+\nu_1(\dD_kG^{ij}_4|\La^{-1}\p_i\p_ja_k)+\nu_1(\La \dD_k G_3|e_k)\nn\\
&&+\nu_1(\dD_kG_4|\La O_k)+2\nu_1( a_k|\p_j \dD_k G^j_0),
\eea
where
\bea
f^2_{\ell,k}&\triangleq&\|a_k\|^2_{L^2}+\|O_k\|^2_{L^2}+\|e_k\|^2_{L^2}+2\nu_1(\La O_k|e_k)-2\nu_1(\La a_k|d_k),\nn\\
\t{f}^2_{\ell,k}&\triangleq&(\mu_0-\nu_1)\|\La e_k\|^2_{L^2}+(\la_0+\mu_0-\nu_1)\|\La d_k\|^2_{L^2}+3\nu_1\|\La a_k\|^2_{L^2}\nn\\
             &&+\nu_1\|\La O_k\|^2_{L^2}+\nu_1\mu_0(\La^2e_k|\La O_k)-\nu_1(\la_0+\mu_0)(\La O^{ij}_k|\p_i\p_j d_k)\nn\\
             &&-\nu_1(\La^2 d_k|\La a_k).\nn
\eea
For any fixed $R_0$, we choose  $\nu_1\sim\nu_1(\la_0,\mu_0,R_0)$ sufficiently small such that
\be\label{e4.11}
\begin{aligned}
&f^2_{\ell,k}\sim \|a_k\|^2_{L^2}+\|e_k\|^2_{L^2}+\|O_k\|^2_{L^2},\\
&\t{f}^2_{\ell,k}\sim 2^{2k}\big(\|a_k\|^2_{L^2}+\|e_k\|^2_{L^2}+\|O_k\|^2_{L^2}\big).
\end{aligned}
\ee
By using Cauchy-Schwarz inequality in \eqref{e4.10}, we can get the following equality owing to $2^k\leq R_0$,
\be\label{e4.12}
\f{d}{dt}f_{\ell,k}+2^{2k}f_{\ell,k}\ls \sum\limits_{i=0,1,3,4}\|\dD_k G_i\|_{L^2},
\ee
which indicates that
\bea\label{e4.13}
\hspace{-15mm}&&\|(a,O;e)\|^{\ell}_{\t{L}^\i_T\dB^{\frac{n}{2}-1}_{2,1}}+\|(a,O;e)\|^{\ell}_{L^1_T\dB^{\frac{n}{2}+1}_{2,1}}
\ls \|(a_0,O_0;e_0)\|^{\ell}_{\dB^{\frac{n}{2}-1}_{2,1}}+\sum\limits_{i=0,1,3,4}\|G_i\|^\ell_{L^1_T\dB^{\frac{n}{2}-1}_{2,1}}.
\eea

Next we begin to bound those nonlinear terms arising in $G_i(i=0,1,3,4)$. Since the quadratic terms containing $a$ and $v$ have already been done in \cite{CMZ1}, it suffices to deal with different terms involving in $O$ as well as those cubic terms due to density-dependent viscosities. More precisely,
we need to estimate the following terms according to the definitions of $G_i$,
\be
G^j_0:=\p_i(aO^{ij}),\ G^{ij}_3:=\p_k v^iO^{kj}-v\cdot\na O^{ij},\nn
\ee
\be\label{4.19er}
\La^{-1}\p_j(O^{lk}\p_lO^{ik}),\ \La^{-1}\p_k(O^{lj}\p_lO^{ik}),\ \La^{-1}\p_k(O^{lk}\p_lO^{ij}) \ \text{in}\ G^{ij}_4,
\ee
and
\be\label{4.20er}
\La^{-1}\p_j\Big(\f{1}{1+a}\text{div}\big(2\t{\mu}(a)D(v)+\t{\la}(a)\text{div} v\text{Id}\big)\Big)^i  \ \text{in}\ G^{ij}_4.
\ee
We write $G^j_0=\p_i aO^{ij}+a\p_iO^{ij}$. Regarding $\p_i aO^{ij}$, by taking $\gamma=-1,r_1=\i,r_2=1,r_3=r_4=2,s_1=s_2=n/2-1,t_1=t_2=n/2$ in \eqref{A.2}
and using \eqref{4.1er}, we arrive at
\bea\label{e5.7}
&&\sum\limits_{{2^k\leq R_0}}2^{k(n/2-1)}\|\dD_k\big(\p_i aO^{ij}\big)\|_{L^1_TL^2}\nn\\
&\ls&\|O\|_{\t{L}^{\i}_T\B^{n/2-1,n/p-1}_{2,p}}\|\na a\|_{\t{L}^{1}_T\B^{n/2,n/p-1}_{2,p}}+\|\na a\|_{\t{L}^{2}_T\B^{n/2-1,n/p-1}_{2,p}}\|O\|_{\t{L}^{2}_T\B^{n/2,n/p}_{2,p}}\nn\\
&\ls&\|O\|_{\t{L}^{\i}_T\B^{n/2-1,n/p}_{2,p}}\| a\|_{\t{L}^{1}_T\B^{n/2+1,n/p}_{2,p}}+\| a\|_{\t{L}^{2}_T\B^{n/2,n/p}_{2,p}}\|O\|_{\t{L}^{2}_T\B^{n/2,n/p}_{2,p}}\nn\\
&\ls&\|(a,O;v)\|^2_{\E^{n/p}_T}.
\eea
The terms $a\p_iO^{ij}$, $v\cdot\na O^{ij}$ in $G^{ij}_3$ and \eqref{4.19er} may be treated along the same lines as $\p_i aO^{ij}$, so we
omit the details for brevity. In order to bound $\p_k v^iO^{kj}$ in $G^{ij}_3$, by taking $\gamma=0,r_1=\i,r_2=1,r_3=r_4=2,s_1=s_2=n/2-1,t_1=t_2=n/2$
in \eqref{A.2} and using \eqref{4.1er}, we have
\bea\label{4.22eer}
&&\sum\limits_{{2^k\leq R_0}}2^{k(n/2-1)}\|\dD_k(\p_k v^iO^{kj})\|_{L^1_TL^2}\nn\\
&\ls&\|O\|_{\t{L}^{\i}_T\B^{n/2-1,n/p-1}_{2,p}}\|\na v\|_{\t{L}^{1}_T\B^{n/2,n/p}_{2,p}}+\|\na v\|_{\t{L}^{2}_T\B^{n/2-1,n/p-1}_{2,p}}\|O\|_{\t{L}^{2}_T\B^{n/2,n/p}_{2,p}}\nn\\
&\ls&\|O\|_{\t{L}^{\i}_T\B^{n/2-1,n/p}_{2,p}}\| v\|_{\t{L}^{1}_T\B^{n/2+1,n/p+1}_{2,p}}+\| v\|_{\t{L}^{2}_T\B^{n/2,n/p}_{2,p}}\|O\|_{\t{L}^{2}_T\B^{n/2,n/p}_{2,p}}\nn\\
&\ls&\|(a,O;v)\|^2_{\E^{n/p}_T}.
\eea

Next we bound the cubic term \eqref{4.20er} in $G^{ij}_4$. Denote
\bea
&&I:=\f{1}{1+a}\text{div}\big(2\t{\mu}(a)D(v)\big)\nn\\
&&\q =\f{1}{1+a}\t{\mu}(a)\na^2 v+\f{1}{1+a}\na\t{\mu}(a)\na v\nn\\
&&\q :=I_1+I_2.\nn
\eea
To bound $I_1$, we have
\bea\label{4.23er}
&&\sum\limits_{{2^k\leq R_0}}2^{k(n/2-1)}\|\dD_k(\f{1}{1+a}\t{\mu}(a)\na^2 v)\|_{L^1_TL^2}\nn\\
&\ls&\sum\limits_{{2^k\leq R_0}}2^{k(n/2-1)}\Big(\|\dD_k(I(a)\t{\mu}(a)\na^2 v)\|_{L^1_TL^2}+\|\dD_k(\t{\mu}(a)\na^2 v)\|_{L^1_TL^2}\Big)\nn\\
&\ls&\|I(a)\|_{\t{L}^{\i}_T\B^{n/2-1,n/p-1}_{2,p}}\|\t{\mu}(a)\na^2 v\|_{\t{L}^{1}_T\B^{n/2,n/p-1}_{2,p}}\nn\\
&&+\|\t{\mu}(a)\na^2 v \|_{\t{L}^{1}_T\B^{n/2-1,n/p-1}_{2,p}}\|I(a)\|_{\t{L}^{\i}_T\B^{n/2,n/p}_{2,p}}+\|\t{\mu}(a)\na^2 v\|_{\t{L}^{1}_T\B^{n/2-1,n/p-1}_{2,p}}\nn\\
&\ls& \Big(1+\|I(a)\|_{\t{L}^{\i}_T\B^{n/2-1,n/p}_{2,p}}\Big)\|\t{\mu}(a)\na^2 v\|_{\t{L}^{1}_T\B^{n/2-1,n/p-1}_{2,p}},
\eea
where we have applied $s_1=s_2=n/2-1, t_1=t_2=n/2, r_1=r_4=\i, r_2=r_3=1,\g=-1$  in \eqref{A.2} of Proposition\ref{A.1} to deal with the term $I(a)\t{\mu}(a)\na^2 v$.
Now we need to estimate $\|\t{\mu}(a)\na^2 v\|_{\t{L}^{1}_T\B^{n/2-1,n/p-1}_{2,p}}$. From \eqref{A.2} and \eqref{A.1}, we have
\bea\label{4.24er}
&&\sum\limits_{{2^k\leq R_0}}2^{k(n/2-1)}\|\dD_k(\t{\mu}(a)\na^2 v)\|_{L^1_TL^2}\nn\\
&\ls&\|\t{\mu}(a)\|_{\t{L}^{\i}_T\B^{n/2-1,n/p-1}_{2,p}}\|\na^2 v\|_{\t{L}^{1}_T\B^{n/2,n/p-1}_{2,p}}+\|\na^2 v \|_{\t{L}^{1}_T\B^{n/2-1,n/p-1}_{2,p}}\|\t{\mu}(a)\|_{\t{L}^{\i}_T\B^{n/2,n/p}_{2,p}}\nn\\
&\ls&\|\t{\mu}(a)\|_{\t{L}^{\i}_T\B^{n/2-1,n/p}_{2,p}}\|v\|_{\t{L}^{1}_T\B^{n/2+1,n/p+1}_{2,p}}.
\eea
\bea\label{4.25er}
&&\sum\limits_{{2^k> R_0}}2^{k(n/p-1)}\|\dD_k(\t{\mu}(a)\na^2 v)\|_{L^1_TL^p}\nn\\
&\ls&\|\t{\mu}(a)\|_{\t{L}^{\i}_T\B^{n/2,n/p}_{2,p}}\|\na^2 v\|_{\t{L}^{1}_T\B^{n/2-1,n/p-1}_{2,p}}.
\eea
Inserting \eqref{4.24er} and \eqref{4.25er} into \eqref{4.23er} and applying Proposition\ref{pa.2}, we can get
\bea\label{4.26er}
&&\sum\limits_{{2^k\leq R_0}}2^{k(n/2-1)}\|\dD_k(\f{1}{1+a}\t{\mu}(a)\na^2 v)\|_{L^1_TL^2}\nn\\
&\ls&(1+\|I(a)\|_{\t{L}^{\i}_T\B^{n/2-1,n/p}_{2,p}})\|\t{\mu}(a)\|_{\t{L}^{\i}_T\B^{n/2-1,n/p}_{2,p}}\| v\|_{\t{L}^{1}_T\B^{n/2+1,n/p+1}_{2,p}}.\nn\\
&\ls&(1+\|a\|_{\t{L}^\i_T\B^{n/p,n/p}_{2,p}})^{n+3}\|a\|_{\t{L}^{\i}_T\B^{n/2-1,n/p}_{2,p}}\|v\|_{\t{L}^{1}_T\B^{n/2+1,n/p+1}_{2,p}}\nn\\
&\ls&(1+\|(a,O;v)\|_{\E^{n/p}_T})^{n+3}\|(a,O;v)\|^2_{\E^{n/p}_T}.
\eea
To bound $I_2$, we have
\bea\label{x.3}
&&\sum\limits_{{2^k\leq R_0}}2^{k(n/2-1)}\|\dD_k(\f{1}{1+a}\na\t{\mu}(a)\na v)\|_{L^1_TL^2}\nn\\
&\ls&\sum\limits_{{2^k\leq R_0}}2^{k(n/2-1)}\Big(\|\dD_k(I(a)\na\t{\mu}(a)\na v)\|_{L^1_TL^2}+\|\dD_k(\na\t{\mu}(a)\na v)\|_{L^1_TL^2}\Big)\nn\\
&\ls&\|I(a)\|_{\t{L}^{\i}_T\B^{n/2-1,n/p-1}_{2,p}}\|\na\t{\mu}(a)\na v\|_{\t{L}^{1}_T\B^{n/2,n/p-1}_{2,p}}\nn\\
&&+\|\na\t{\mu}(a)\na v \|_{\t{L}^{1}_T\B^{n/2-1,n/p-1}_{2,p}}\|I(a)\|_{\t{L}^{\i}_T\B^{n/2,n/p}_{2,p}}+\|\na\t{\mu}(a)\na v\|_{\t{L}^{1}_T\B^{n/2-1,n/p-1}_{2,p}}\nn\\
&\ls& \Big(1+\|I(a)\|_{\t{L}^{\i}_T\B^{n/2-1,n/p}_{2,p}}\Big)\|\na\t{\mu}(a)\na v\|_{\t{L}^{1}_T\B^{n/2-1,n/p-1}_{2,p}},
\eea
 From \eqref{A.2} and \eqref{A.1}, the estimate of $\|\na\t{\mu}(a)\na v\|_{\t{L}^{1}_T\B^{n/2-1,n/p-1}_{2,p}}$ is
\bea\label{x.1}
&&\sum\limits_{{2^k\leq R_0}}2^{k(n/2-1)}\|\dD_k(\na\t{\mu}(a)\na v)\|_{L^1_TL^2}\nn\\
&\ls&\|\na v\|_{\t{L}^{2}_T\B^{n/2-1,n/p-1}_{2,p}}\|\na\t{\mu}(a) \|_{\t{L}^{2}_T\B^{n/2,n/p-1}_{2,p}}+\|\na \t{\mu}(a) \|_{\t{L}^{\i}_T\B^{n/2-1,n/p-1}_{2,p}}\|\na v\|_{{L}^{1}_T\B^{n/2,n/p}_{2,p}}\nn\\
&\ls&\|v\|_{\t{L}^{2}_T\B^{n/2,n/p}_{2,p}}\|\t{\mu}(a) \|_{\t{L}^{2}_T\B^{n/2,n/p}_{2,p}}+\| \t{\mu}(a) \|_{\t{L}^{\i}_T\B^{n/2-1,n/p}_{2,p}}\|v\|_{{L}^{1}_T\B^{n/2+1,n/p+1}_{2,p}}.
\eea
\bea\label{x.2}
&&\sum\limits_{{2^k> R_0}}2^{k(n/p-1)}\|\dD_k(\na \t{\mu}(a)\na v)\|_{L^1_TL^p}\nn\\
&\ls&\|\na\t{\mu}(a)\|_{\t{L}^{2}_T\B^{n/2,n/p}_{2,p}}\|\na v\|_{\t{L}^{2}_T\B^{n/2-1,n/p-1}_{2,p}}\nn\\
&\ls&\|v\|_{\t{L}^{2}_T\B^{n/2,n/p}_{2,p}}\|\t{\mu}(a) \|_{\t{L}^{2}_T\B^{n/2,n/p}_{2,p}}
\eea
Inserting \eqref{x.1} and \eqref{x.2} into \eqref{x.3} and applying Proposition\ref{pa.2} and \eqref{4.1er}, we can get
\bea
&&\sum\limits_{{2^k\leq R_0}}2^{k(n/2-1)}\|\dD_k(\f{1}{1+a}\na\t{\mu}(a)\na v)\|_{L^1_TL^2}\nn\\
&\ls&(1+\|I(a)\|_{\t{L}^{\i}_T\B^{n/2-1,n/p}_{2,p}})\Big(\|\t{\mu}(a)\|_{\t{L}^{\i}_T\B^{n/2-1,n/p}_{2,p}}\| v\|_{\t{L}^{1}_T\B^{n/2+1,n/p+1}_{2,p}}.\nn\\
&&\qq+\|v\|_{\t{L}^{2}_T\B^{n/2,n/p}_{2,p}}\|\t{\mu}(a) \|_{\t{L}^{2}_T\B^{n/2,n/p}_{2,p}}\Big)\nn\\
&\ls&(1+\|(a,O;v)\|_{\E^{n/p}_T})^{n+3}\|(a,O;v)\|^2_{\E^{n/p}_T}.
\eea

Since bound of the  cubic term $\f{1}{1+a}\text{div}\big(\t{\la}(a)\text{div} v\text{Id}\big)$ is the same as $I$, we omit the details. Summing up all the estimates and remembering \eqref{e4.13}, we get
\bea\label{4.27er}
&&\|(a,O;e)\|^{\ell}_{\t{L}^\i_T\dB^{n/2-1}_{2,1}}+\|(a,O;e)\|^{\ell}_{L^1_T\dB^{n/2+1}_{2,1}}\nn\\
&\ls& \|(a_0,O_0;e_0)\|^{\ell}_{\dB^{n/2-1}_{2,1}}+(1+\|(a,O;v)\|_{\E^{n/p}_T})^{n+3}\|(a,O;v)\|^2_{\E^{n/p}_T}.
\eea
\noindent\textbf{Step 2: High-frequency estimates ($2^k> R_0$).}

\noindent Inspired by \cite{Hb1,Hb2}, we perform basic energy approaches in terms of \textit{effective
velocities} rather than the Lagrangian change as in \cite{CD1,CMZ1}. Denote by  $\t{d}=-\na (-\Dl)^{-1}\text{div} v$
the compressible part of $v$. It is easy to see that $\|\t{d}\|_{\t{L}^q_T\B^{s,\sigma}_{2,p}}\thickapprox\|d\|_{\t{L}^q_T\B^{s,\sigma}_{2,p}}$.
It follows from the first equality in \eqref{ip1} that
\bea\label{3.8e}
&&-\na(-\Dl)^{-1}\text{div}(\na\cdot O)\nn\\
&=&-\na(-\Dl)^{-1}\Big(\p_i\p_j[(1+a)(\dl^{ij}+O^{ij})]\Big)+\na(-\Dl)^{-1}\text{div}\text{div}(aI+aO)\nn\\
&=&\na(-\Dl)^{-1}\text{div}\text{div}(aI+aO)\nn\\
&=&-\na a+\na(-\Dl)^{-1}\text{div}\text{div}(aO).
\eea
Note that \eqref{3.8e}, we get the following equation for the compressible part of $v$
\be\label{3.9e}
\p_t \t{d}-\Dl \t{d}+2\na a=G_2,
\ee
where
\bea\label{3.10e}
&&G_2=-\na(-\Dl)^{-1}\text{div}\Big(-v\cdot\nabla v+O^{jk}\p_jO^{\bullet k}-I(a)\mathcal{A}v\nn\\
&&\qq-K(a)\nabla a -\text{div}(aO)+\f{1}{1+a}\text{div}\big(2\t{\mu}(a)D(v)+\t{\la}(a)\text{div} v\text{Id}\big)\Big).
\eea
The motivation using the system \eqref{3.9e} is to make a comparison with the usual compressible Navier-Stokes equations. Here, we consider more complicated hyperbolic-parabolic coupled system
\be\label{e4.17}
\left\{
\begin{aligned}
&\p_t a+v\cdot\na a+\text{div}v=\t{G}_1, \\
&\p_t \t{d}-\Dl \t{d}+2\na a=G_2,\\
&\p_tO^{ij}+v\cdot \na O^{ij}-\La e^{ij}=\t{G}^{ij}_3,\\
&\p_t e^{ij}-\mu_0\Dl e^{ij}+\La O^{ij}=\t{G}^{ij}_4\\
\end{aligned}
\right.
\ee
where
\be\label{e4.18}
\t{G}_1=-a\na\cdot v,\q \t{G}^{ij}_3=\p_kv^iO^{kj},\q \t{G}^{ij}_{4}=G^{ij}_{4}+(\la_0+\mu_0)\p_i\p_j d-\La^{-1}\p_i\p_j a.
\ee
Introduce two \textit{effective velocities} as follows
\be
w=\t{d}+2\na(-\Dl)^{-1}a=\na(-\Dl)^{-1}(2a-\text{div} v),\qq \O^{ij}= e^{ij}+\f{1}{\mu_0}\La(-\Dl)^{-1} O^{ij}.\nn
\ee
Noticing that the definition of $w$ is almost the same as that in \cite{Hb1,Hb2}. The subtle difference lies on the coefficient of unknown $a$, which comes from the contribution of deformation gradient $F$, see (\ref{3.8e}). The new \textit{effective velocity} $\O^{ij}$ is used to cancel the coupling between $e^{ij}$ and $O^{ij}$ in the high-frequency estimate.

Firstly, we present those estimates for effective velocities. It follows from \eqref{e4.17} that
\be\label{e4.19}
\left\{
\begin{aligned}
&\p_t w-\Dl w=G_2+2\na(-\Dl)^{-1}\t{G}_1+2w-4\na(-\Dl)^{-1}a, \\
&\p_t \O^{ij}-\mu_0\Dl \O^{ij}=\t{G}^{ij}_4+\f{1}{\mu_0}\La^{-1}\t{G}^{ij}_3+\f{1}{\mu_0}\O^{ij}-\f{1}{\mu^2_0}\La^{-1}O^{ij}.
\end{aligned}
\right.
\ee
Applying \eqref{A.6} to the above equations implies that
\be\label{ee4.20}
\begin{aligned}
&\|w\|^h_{\t{L}^\i_T\dB^{n/p-1}_{p,1}\cap L^1_T\dB^{n/p+1}_{p,1}}\ls \|w_0\|^h_{\dB^{n/p-1}_{p,1}}+\|w\|^h_{L^1_T\dB^{n/p-1}_{p,1}}+\| a\|^h_{L^1_T\dB^{n/p-2}_{p,1}}\\
&\q\qq\qq\qq\qq\hspace{5mm}+\|\t{G}_1\|^h_{L^1(\dB^{n/p-2}_{p,1})}+\|G_2\|^h_{L^1(\dB^{n/p-1}_{p,1})},\\
\end{aligned}
\ee
and
\be\label{ee4.21}
\begin{aligned}
&\|\O^{ij}\|^h_{\t{L}^\i_T\dB^{n/p-1}_{p,1}\cap L^1_T\dB^{n/p+1}_{p,1}}\ls \|\O^{ij}_0\|^h_{\dB^{n/p-1}_{p,1}}+\|\O^{ij}\|^h_{L^1_T\dB^{n/p-1}_{p,1}}+\| O^{ij}\|^h_{L^1_T\dB^{n/p-2}_{p,1}}\\
&\q\qq\qq\qq\qq\hspace{5mm}+\|\t{G}^{ij}_3\|^h_{L^1(\dB^{n/p-2}_{p,1})}+\|\t{G}^{ij}_4\|^h_{L^1(\dB^{n/p-1}_{p,1})}.\\
\end{aligned}
\ee
Owing to the high frequency cut-off $2^k>R_0$, we have
\be
\|w\|^h_{L^1_T\dB^{n/p-1}_{p,1}}\ls R^{-2}_0\|w\|^h_{L^1_T\dB^{n/p+1}_{p,1}}, \q \| a\|^h_{L^1_T\dB^{n/p-2}_{p,1}}\ls R^{-2}_0\|a\|^h_{L^1_T\dB^{n/p}_{p,1}},\nn
\ee
and
\be
\|\O^{ij}\|^h_{L^1_T\dB^{n/p-1}_{p,1}}\ls R^{-2}_0\|w^{ij}\|^h_{L^1_T\dB^{n/p+1}_{p,1}}, \q \| O^{ij}\|^h_{L^1_T\dB^{n/p-2}_{p,1}}\ls R^{-2}_0\|O^{ij}\|^h_{L^1_T\dB^{n/p}_{p,1}}.\nn
\ee
Choosing $R_0>0$ sufficient large, the terms $\|w\|^h_{L^1_T\dB^{n/p-1}_{p,1}}$ and $\|\O^{ij}\|^h_{L^1_T\dB^{n/p-1}_{p,1}}$ on the right-side of
\eqref{ee4.20} and \eqref{ee4.21} can be absorbed by the corresponding parts in the left-hand side. Consequently, we conclude that
\be\label{ee4.22}
\begin{aligned}
&\|w\|^h_{\t{L}^\i_T\dB^{n/p-1}_{p,1}\cap L^1_T\dB^{n/p+1}_{p,1}}\ls \|w_0\|^h_{\dB^{n/p-1}_{p,1}}+R^{-2}_0\| a\|^h_{L^1_T\dB^{n/p}_{p,1}}\\
&\q\qq\qq\qq\qq+\|\t{G}_1\|^h_{L^1(\dB^{n/p-2}_{p,1})}+\|G_2\|^h_{L^1(\dB^{n/p-1}_{p,1})},\\
\end{aligned}
\ee
and
\be\label{ee4.23}
\begin{aligned}
&\|\O^{ij}\|^h_{\t{L}^\i_T\dB^{n/p-1}_{p,1}\cap L^1_T\dB^{n/p+1}_{p,1}}\ls \|\O^{ij}_0\|^h_{\dB^{n/p-1}_{p,1}}+R^{-2}_0\| O^{ij}\|^h_{L^1_T\dB^{n/p}_{p,1}}\\
&\q\qq\qq\qq\qq+\|\t{G}^{ij}_3\|^h_{L^1(\dB^{n/p-2}_{p,1})}+\|\t{G}^{ij}_4\|^h_{L^1(\dB^{n/p-1}_{p,1})}.\\
\end{aligned}
\ee

Secondly, we see that $(a, O^{ij})$ satisfies the following damped equations in terms of effective velocities
\be
\lt\{
\begin{aligned}
&\p_t a+v\cdot\na a+2a=\t{G}_1-\na\cdot w,\\
&\p_t O^{ij}+v\cdot\na O^{ij}+\f{1}{\mu_0}O^{ij}=\t{G}^{ij}_3+\La \O^{ij}.
\end{aligned}
\rt.\nn
\ee
Applying $\dD_k$ to the above equations, we obtain
\be\label{ee4.24}
\lt\{
\begin{aligned}
&\p_t \dD_ka+v\cdot\na \dD_ka+2\dD_ka=\dD_k\t{G}_1-\dD_k\na\cdot w+R^1_k,\\
&\p_t \dD_kO^{ij}+v\cdot\na \dD_kO^{ij}+\f{1}{\mu_0}\dD_kO^{ij}=\dD_k\t{G}^{ij}_3+\dD_k\La \O^{ij}+R^2_k,
\end{aligned}
\rt.
\ee
where $R^1_k:=[v\cdot\na, \dD_k]a$ and $R^2_k:=[v\cdot\na, \dD_k]O^{ij}$.
Multiplying $\eqref{ee4.24}_1$ by $\dD_ka|\dD_ja|^{p-2}$ and $\eqref{ee4.24}_2$ by $\dD_kO^{ij}|\dD_kO^{ij}|^{p-2}$, and then integrating over $\bR^n\times [0,t]$, we can obtain
\bea\label{ee4.25}
&&\|\dD_ka(t)\|_{L^p}+\int^t_0\|\dD_ka\|_{L^p}d\tau\ls \|\dD_k a_0\|_{L^p}+\int^t_0\|\na v\|_{L^\i}\|\dD_ka\|_{L^p}d\tau\nn\\
&&\qq+\int^t_0\|\dD_k(\t{G}_1-\La w)\|_{L^p}d\tau+\int^t_0\|R^1_k\|_{L^p}d\tau
\eea
and
\bea\label{ee4.26}
&&\|\dD_kO^{ij}(t)\|_{L^p}+\int^t_0\|\dD_kO^{ij}\|_{L^p}d\tau\ls \|\dD_k O^{ij}_0\|_{L^p}+\int^t_0\|\na v\|_{L^\i}\|\dD_kO^{ij}\|_{L^p}d\tau\nn\\
&&\qq+\int^t_0\|\dD_k(\t{G}^{ij}_3+\La \O^{ij})\|_{L^p}d\tau+\int^t_0\|R^2_k\|_{L^p}d\tau.
\eea
It follows from commutator estimates in \cite{BCD1} that
\be
\sum\limits_{j\in\bZ}2^{js}\|(R^1_j,R^2_j)\|_{L^p}\ls \|\na v\|_{\dB^{n/p}_{p,1}}\|(a,O^{ij})\|_{\dB^{s}_{p,1}}.\nn
\ee
Now multiplying \eqref{ee4.25} and \eqref{ee4.26} by $2^{k\f{n}{p}}$, respectively,  and then summing over the index $k$ satisfying
$2^k>R_0$, we are led to
\bea\label{ee4.27}
\|a\|^h_{\t{L}^\i_T\dB^{n/p}_{p,1}\cap{L}^1_T\dB^{n/p}_{p,1}}&\ls& \|a_0\|^h_{\dB^{n/p}_{p,1}}+\|\na v\|_{L^1_T\dB^{n/p}_{p,1}}\|a\|_{\t{L}^\i_T\dB^{n/p}_{p,1}}+\|\t{G}_1\|^h_{L^1_T\dB^{n/p}_{p,1}}+\|w\|^h_{L^1_T(\dB^{n/p+1}_{p,1})}\nn\\
&\ls&\|a_0\|^h_{\dB^{n/p}_{p,1}}+\|w\|^h_{L^1_T(\dB^{n/p+1}_{p,1})}+\|(a,O;v)\|^2_{\E^{n/p}_T}+\|\t{G}_1\|^h_{L^1_T\dB^{n/p}_{p,1}}
\eea
and
\bea\label{ee4.28}
&&\|O^{ij}\|^h_{\t{L}^\i_T\dB^{s}_{p,1}\cap{L}^1_T\dB^{n/p}_{p,1}}\nonumber\\ &\ls& \|O^{ij}_0\|^h_{\dB^{n/p}_{p,1}}+\|\na v\|_{L^1_T\dB^{n/p}_{p,1}}\|O^{ij}\|_{\t{L}^\i_T\dB^{n/p}_{p,1}}+\|\t{G}^{ij}_3\|^h_{L^1_T\dB^{n/p}_{p,1}}+\|\O^{ij}\|^h_{L^1_T(\dB^{n/p+1}_{p,1})}\nn\\
&\ls&\|O^{ij}_0\|^h_{\dB^{n/p}_{p,1}}+\|\O^{ij}\|^h_{L^1_T(\dB^{n/p+1}_{p,1})}+\|(a,O;v)\|^2_{\E^{n/p}_T}+\|\t{G}^{ij}_3\|^h_{L^1_T\dB^{n/p}_{p,1}}.
\eea
Multiply \eqref{ee4.27} and \eqref{ee4.28} by $\dl>0$ respectively, and then add two resulting inequalities to \eqref{ee4.22} and \eqref{ee4.23} together. By choosing $R_0$ sufficiently large, we get
 \bea
&&\|a\|^h_{\t{L}^\i_T\dB^{n/p}_{p,1}\cap{L}^1_T\dB^{n/p}_{p,1}}+\|w\|^h_{\t{L}^\i_T\dB^{n/p-1}_{p,1}\cap L^1_T\dB^{n/p+1}_{p,1}}\ls\|a_0\|^h_{\dB^{n/p}_{p,1}}+\|w_0\|^h_{\dB^{n/p-1}_{p,1}}\nn\\
&&\qq\qq\qq +\|(a,O;v)\|^2_{\E^{n/p}_T}+\|\t{G}_1\|^h_{L^1_T(\dB^{n/p}_{p,1})}+\|G_2\|^h_{L^1_T(\dB^{n/p-1}_{p,1})},\nn
 \eea
 and
 \bea
&&\|O^{ij}\|^h_{\t{L}^\i_T\dB^{n/p}_{p,1}\cap{L}^1_T\dB^{n/p}_{p,1}}+\|\O^{ij}\|^h_{\t{L}^\i_T\dB^{n/p-1}_{p,1}\cap L^1_T\dB^{n/p+1}_{p,1}}\ls\|O^{ij}_0\|^h_{\dB^{n/p}_{p,1}}+\|\O^{ij}_0\|^h_{\dB^{n/p-1}_{p,1}}\nn\\
&&\qq\qq\qq +\|(a,O;v)\|^2_{\E^{n/p}_T}+\|\t{G}^{ij}_3\|^h_{L^1_T(\dB^{n/p}_{p,1})}+\|\t{G}^{ij}_4\|^h_{L^1_T(\dB^{n/p-1}_{p,1})}.\nn
 \eea
Keep in mind that $w=\t{d}+2\na(-\Dl)^{-1}a,\ \O^{ij}=e^{ij}+\f{1}{\mu_0}(-\Dl)^{-1}\La O^{ij}$, we arrive at
\bea\label{ee4.29}
&&\|a\|^h_{\t{L}^\i_T\dB^{n/p}_{p,1}\cap{L}^1_T\dB^{n/p}_{p,1}}+\|\t{d}\|^h_{\t{L}^\i_T\dB^{n/p-1}_{p,1}\cap L^1_T\dB^{n/p+1}_{p,1}}\ls\|a_0\|^h_{\dB^{n/p}_{p,1}}+\|\t{d}_0\|^h_{\dB^{n/p-1}_{p,1}}\nn\\
&&\qq\qq\qq +\|(a,O;v)\|^2_{\E^{n/p}_T}+\|\t{G}_1\|^h_{L^1_T(\dB^{n/p}_{p,1})}+\|G_2\|^h_{L^1_T(\dB^{n/p-1}_{p,1})}
 \eea
 and
 \bea\label{ee4.30}
&&\|O^{ij}\|^h_{\t{L}^\i_T\dB^{n/p}_{p,1}\cap{L}^1_T\dB^{n/p}_{p,1}}+\|e^{ij}\|^h_{\t{L}^\i_T\dB^{n/p-1}_{p,1}\cap L^1_T\dB^{n/p+1}_{p,1}}\ls\|O^{ij}_0\|^h_{\dB^{n/p}_{p,1}}+\|e^{ij}_0\|^h_{\dB^{n/p-1}_{p,1}}\nn\\
&&\qq\qq\qq +\|(a,O;v)\|^2_{\E^{n/p}_T}+\|\t{G}^{ij}_3\|^h_{L^1_T(\dB^{n/p}_{p,1})}+\|\t{G}^{ij}_4\|^h_{L^1_T(\dB^{n/p-1}_{p,1})}.
 \eea
In addition, remembering \eqref{e4.18}, we have
\be\label{ee4.31}
\|\t{G}^{ij}_4\|^h_{L^1_T(\dB^{n/p-1}_{p,1})}\ls \|{G}^{ij}_4\|^h_{L^1_T(\dB^{n/p-1}_{p,1})}+\|\t{d}\|^h_{L^1_T(\dB^{n/p+1}_{p,1})}+\|a\|^h_{L^1(\dB^{n/p}_{p,1})}
\ee
Hence, together with \eqref{ee4.29}-\eqref{ee4.31}, we deduce that
 \bea\label{ee4.32}
&&\|(a,O)\|^h_{\t{L}^\i_T\dB^{n/p}_{p,1}\cap{L}^1_T\dB^{n/p}_{p,1}}+\|e\|^h_{\t{L}^\i_T\dB^{n/p-1}_{p,1}\cap L^1_T\dB^{n/p+1}_{p,1}}\nn\\
&\ls&\|(a_0,O_0)\|^h_{\dB^{n/p}_{p,1}}+\|e_0\|^h_{\dB^{n/p-1}_{p,1}}+\|(a,O;v)\|^2_{\E^{n/p}_T}+\|(\t{G}_1,\t{G}_3)\|^h_{L^1_T(\dB^{n/p}_{p,1})}\nn\\
&&+\|(G_2,G_4)\|^h_{L^1_T(\dB^{n/p-1}_{p,1})}.
 \eea

Likely, we need to bound those different terms in $\t{G}_i(i=1,3)$ and $G_i(i=2,4)$ compared to \cite{CMZ1}, for example,
\be
 \t{G}^{ij}_3:=\p_k v^iO^{kj},\nn
\ee
\be
O\na O,\ \text{div}(aO) \ \text{in}\ G_2\ \text{and}\ G_4,\nn
\ee
and
\be\label{4.47er}
\f{1}{1+a}\text{div}\big(2\t{\mu}(a)D(v)+\t{\la}(a)\text{div} v\text{Id}\big)  \ \text{in}\ G_2\ \text{and}\ G_4.
\ee

In order to bound $\p_k v^iO^{kj}$, from \eqref{A.1} of Proposition\ref{A.1} with $,r_1=1,r_2=\i,\sigma=\tau=n/p$, we have

\bea
&&\sum\limits_{{2^k> R_0}}2^{k(n/p)}\|\dD_k(\p_k v^iO^{kj})\|_{L^1_TL^p}\nn\\
&\ls&\|\na v\|_{\t{L}^{1}_T\B^{n/2,n/p}_{2,p}}\|O\|_{\t{L}^{\i}_T\B^{n/2,n/p}_{2,p}}\ls\|v\|_{\t{L}^{1}_T\B^{n/2+1,n/p+1}_{2,p}}\|O\|_{\t{L}^{\i}_T\B^{n/2-1,n/p}_{2,p}}\nn\\
&\ls&\|(a,O;v)\|^2_{\E^{n/p}_T}.\nn
\eea

For $O\na O$,  from \eqref{A.1} of Proposition\ref{A.1} with $r_1=r_2=2,\sigma=n/p,\tau=n/p-1$ and by applying interpolation \eqref{4.1er}, we have

\bea\label{4.48rt}
&&\sum\limits_{{2^k> R_0}}2^{k(n/p-1)}\|\dD_k\big(O\na O\big)\|_{L^1_TL^p}\nn\\
&\ls&\|O\|_{\t{L}^{2}_T\B^{n/2,n/p}_{2,p}}\|\na O\|_{\t{L}^{2}_T\B^{n/2-1,n/p-1}_{2,p}}\ls\|O\|^2_{\t{L}^{2}_T\B^{n/2,n/p}_{2,p}} \nn\\
&\ls&\|(a,O;v)\|^2_{\E^{n/p}_T}.
\eea

The estimate for $\text{div}(aO)=a\na\cdot O+\na a O$ may be handled with at the same away as $O\na O$. Next, we bound the cubic term \eqref{4.47er} in $G^{ij}_4$. Following from the the same notation, we know that
\bea
&&I=\f{1}{1+a}\text{div}\big(2\t{\mu}(a)D(v)\big)\nn\\
&&\q =\f{1}{1+a}\t{\mu}(a)\na^2 v+\f{1}{1+a}\na\t{\mu}(a)\na v\nn\\
&&\q \triangleq I_1+I_2.\nn
\eea

For $I_1$, we have
\bea
&&\sum\limits_{{2^k> R_0}}2^{k(n/p-1)}\|\dD_k(\f{1}{1+a}\t{\mu}(a)\na^2 v)\|_{L^1_TL^p}\nn\\
&\ls&\sum\limits_{{2^k> R_0}}2^{k(n/p-1)}\Big(\|\dD_k(I(a)\t{\mu}(a)\na^2 v)\|_{L^1_TL^p}+\|\dD_k(\t{\mu}(a)\na^2 v)\|_{L^1_TL^p}\Big)\nn\\
&\ls&\|I(a)\|_{\t{L}^{\i}_T\B^{n/2,n/p}_{2,p}}\|\t{\mu}(a)\na^2 v\|_{\t{L}^{1}_T\B^{n/2-1,n/p-1}_{2,p}}+\|\t{\mu}(a)\na^2 v\|_{\t{L}^{1}_T\B^{n/2-1,n/p-1}_{2,p}}\nn\\
&\ls& \Big(1+\|I(a)\|_{\t{L}^{\i}_T\B^{n/2-1,n/p}_{2,p}}\Big)\|\t{\mu}(a)\na^2 v\|_{\t{L}^{1}_T\B^{n/2-1,n/p-1}_{2,p}}\nn\\
&\ls&(1+\|(a,O;v)\|_{\E^{n/p}_T})^{n+3}\|(a,O;v)\|^2_{\E^{n/p}_T},\nn
\eea
where the third line is followed by taking $\sigma=n/p,\tau=n/p-1, r_1=\i, r_2=1$  in \eqref{A.1}.

On the other hand, regarding $I_2$, we deduce that
\bea
&&\sum\limits_{{2^k> R_0}}2^{k(n/p-1)}\|\dD_k(\f{1}{1+a}\na\t{\mu}(a)\na v)\|_{L^1_TL^p}\nn\\
&\ls&\sum\limits_{{2^k> R_0}}2^{k(n/p-1)}\Big(\|\dD_k(I(a)\na\t{\mu}(a)\na v)\|_{L^1_TL^p}+\|\dD_k(\na\t{\mu}(a)\na v)\|_{L^1_TL^p}\Big)\nn\\
&\ls&\|I(a)\|_{\t{L}^{\i}_T\B^{n/2,n/p}_{2,p}}\|\na\t{\mu}(a)\na v\|_{\t{L}^{1}_T\B^{n/2-1,n/p-1}_{2,p}}+\|\na\t{\mu}(a)\na v\|_{\t{L}^{1}_T\B^{n/2-1,n/p-1}_{2,p}}\nn\\
&\ls& \Big(1+\|I(a)\|_{\t{L}^{\i}_T\B^{n/2-1,n/p}_{2,p}}\Big)\|\na\t{\mu}(a)\na v\|_{\t{L}^{1}_T\B^{n/2-1,n/p-1}_{2,p}}\nn\\
&\ls&(1+\|(a,O;v)\|_{\E^{n/p}_T})^{n+3}\|(a,O;v)\|^2_{\E^{n/p}_T}.\nn
\eea

The computation for $\f{1}{1+a}\text{div}\big(\t{\la}(a)\text{div} v\text{Id}\big)$ totally follows from the same procedure as $I$, so we omit details. By putting above estimates together, remembering \eqref{ee4.32}, we achieve that
 \bea\label{4.51er}
&&\|(a,O)\|^h_{\t{L}^\i_T\dB^{n/p}_{p,1}\cap{L}^1_T\dB^{n/p}_{p,1}}+\|e\|^h_{\t{L}^\i_T\dB^{n/p-1}_{p,1}\cap L^1_T\dB^{n/p+1}_{p,1}}\nn\\
&\ls&\|(a_0,O_0)\|^h_{\dB^{n/p}_{p,1}}+\|e_0\|^h_{\dB^{n/p-1}_{p,1}}+(1+\|(a,O;v)\|_{\E^{n/p}_T})^{n+3}\|(a,O;v)\|^2_{\E^{n/p}_T}.
 \eea
\noindent\textbf{Step 3: Combination of two-step analysis.}

\noindent The inequality (\ref{4.1}) is the consequence of \eqref{4.27er} and \eqref{4.51er}, so the proof of Proposition \ref{pro4.1} is finished. By using \eqref{A.3} in Proposition \ref{A.1}, we can infer that
\bea\label{4.52er}
&&\|(a,O;e)\|^{\ell}_{\t{L}^\i_T\dB^{n/2-1}_{2,1}}+\|(a,O;e)\|^{\ell}_{L^1_T\dB^{n/2+1}_{2,1}}\nn\\
&\ls& \|(a_0,O_0;e_0)\|^{\ell}_{\dB^{n/2-1}_{2,1}}\nn\\
&&+(1+\|(a,O;v)\|_{\E^{n/p}_T})^{n+3}\|(a,O;v)\|_{\E^{n/p}_T}\|(a,O;v)\|_{E^{n/p}_T}
\eea
and
\bea\label{4.53er}
&&\|(a,O)\|^h_{\t{L}^\i_T\dB^{n/2}_{2,1}\cap{L}^1_T\dB^{n/2}_{2,1}}+\|e\|^h_{\t{L}^\i_T\dB^{n/2-1}_{p,1}\cap L^1_T\dB^{n/2+1}_{p,1}}\nn\\
&\ls&\|(a_0,O_0)\|^h_{\dB^{n/2}_{2,1}}+\|e_0\|^h_{\dB^{n/2-1}_{2,1}}\nn\\
&&+(1+\|(a,O;v)\|_{\E^{n/p}_T})^{n+3}\|(a,O;v)\|_{\E^{n/p}_T}\|(a,O;v)\|_{E^{n/2}_T}.
\eea
The inequality (\ref{4.8r}) is followed by \eqref{4.52er} and \eqref{4.53er}. Therefore, the proof of
Proposition \ref{pro4.2} is complete.

\subsection{Approximate solutions and uniform estimates}
The construction of approximate solutions is based on the following local-in-time existence.
\begin{theorem}[\cite{QZ1}]\label{th4.1}
Assume $(\rho_0-1,F_0-I)\in \big(\dot{B}^{n/2}_{2,1}\big)^{1+n^2}$ and $u_0\in\big(\dot{B}^{n/2-1}_{2,1}\big)^n$ with $\rho_0$ bounded away from $0$. There exists a positive time $T$ such that system \eqref{1.1} has a unique solution $(\rho,F;u)$ with $\rho$ bounded away from $0$ and
\be
(\rho-1,F-I)\in \Big(C([0,T);\dot{B}^{n/2}_{2,1})\Big)^{1+n^2}, u\in \Big(C([0,T);\dot{B}^{n/2-1}_{2,1})\cap L^1([0,T);\dot{B}^{n/2+1}_{2,1})\Big)^n
.\nn
\ee
Additionally, if $(\rho_0-1,F_0-I)\in\big(\dot{B}^{n/2-1}_{2,1}\big)^{1+n^2}$, we have
\be
(\rho-1,F-I)\in\big(C([0,T);\dot{B}^{n/2-1}_{2,1})\big)^{1+n^2}.\nn
\ee
\end{theorem}
In order to apply Theorem \ref{th4.1}, we need the following lemma, which could be shown by the proof of Lemma 4.2 in \cite{Ah1}.
\begin{lemma}\label{l5.1}
Let $p\geq 2$. For any
\be
(\rho_0-1,F_0-I;u_0)\in \Big(\B^{n/2-1,n/p}_{2,p}\Big)^{1+n^2}\times\Big(\B^{n/2-1,n/p-1}_{2,p}\Big)^n\nn
\ee
satisfying $\rho_0\geq c_0>0$, then there exists a sequence $\{(\rho_{0,k},F_{0,k};u_{0,k})\}_{k\in\bN}$ with $\{(\rho_{0,k}-1,F_{0,k}-I;u_{0,k})\}
\in \Big(\B^{n/2-1,n/2}_{2,2}\Big)^{1+n^2}\times\Big(\dot{B}^{n/2-1}_{2,1}\Big)^n $ such that
\be\label{4.14}
\|(\rho_{0,k}-\rho_0,F_{0,k}-F_0)\|_{\B^{n/2-1,n/p}_{2,p}}\longrightarrow 0,\q \|u_{0,k}-u_0\|_{\B^{n/2-1,n/p-1}_{2,p}}\longrightarrow 0
\ee
when $k\rightarrow 0$. we also have $\rho_{0,k}\geq \f{c_0}{2}$ for any $k\in\bN$.
\end{lemma}

Let $(\rho_{0,k},F_{0,k};u_{0,k})$ be the sequence for initial data stated in Lemma \ref{l5.1}. Then Theorem \ref{th4.1} indicates that there exists a maximal existence time $T_k>0$ such that
system \eqref{1.1} with initial data $(\rho_{0,k},F_{0,k};u_{0,k})$ has a unique solution $(\rho_k,F_k;u_k)$ with $\rho_k$ bounded away from $0$, and satisfies
\bea
&&(\rho_k-1,F_k-I)\in \Big(C([0,T_k);\dot{B}^{n/2}_{2,1}\cap\dot{B}^{n/2-1}_{2,1})\Big)^{1+n^2},\nn\\
&&u_k\in\Big(C([0,T_k);\dot{B}^{n/2-1}_{2,1})\cap L^1(\dot{B}^{n/2+1}_{2,1})\Big)^{n}.\nn
\eea
Then using the definition of Hybird Besov spaces and Bernstein inequality in Lemma \ref{lem2.1}, we have
\bea
&&(\rho_k-1,F_k-I)\in \Big(C([0,T_k);\B^{n/2-1,n/p}_{2,p})\Big)^{1+n^2},\nn\\
&&u_k\in\Big(C([0,T_k);\B^{n/2-1,n/p-1}_{2,p})\cap L^1([0,T_k);\B^{n/2+1,n/p+1}_{2,p})\Big)^{n}.\nn
\eea
Set
\be
a_k(t,x)=\rho_k(\chi^2_0t,\chi_0x)-1,\q v_k(t,x)=\chi_0u_k(\chi^2_0t,\chi_0x),\q O_k(t,x)=F_k(\chi^2_0t,\chi_0x)-I.\nn
\ee
From \eqref{1.3rx} and \eqref{4.14}, we
\be
\|(a_{0,k},O_{0,k};v_{0,k})\|_{\mathcal{E}^{n/p}_0}\leq C_0\eta,\nn
\ee
for some constant $C_0$. Given a constant $M$ to be determined later on, we define
\be
T^\ast_k \triangleq \sup\{t\in[0,T_k)\big| \|(a_k,O_k;v_k)\|_{\mathcal{E}^{n/p}_t}\leq M\eta\}.\nn
\ee
First we claim that
\be
T^\ast_k=T_k \qq \forall k\in \bN.\nn
\ee
With aid of the continuity argument, it suffices to show for all $k\in \bN$,
\be\label{4.15}
\|(a_k,O_k;v_k)\|_{\mathcal{E}^{n/p}_{T^\ast_k}}\leq \f{1}{2}M\eta.
\ee
Indeed, noting that $\|a_k\|_{L^\i([0,T^\ast_k)\times\bR^n)}\leq C_1\|a^k\|_{L^\i_{T^\ast_k}{\B^{n/2-1,n/p}_{2,p}}}$, we can choose  $\eta$ sufficiently small such that
\be
M\eta\leq \f{1}{2C_1}.\nn
\ee
Then
\be
\|a_k\|_{L^\i([0,T^\ast_k)\times\bR^n)}\leq \f{1}{2}.\nn
\ee
By applying Proposition \ref{pro4.1}, we obtain
\begin{equation}\label{4.155}
\|(a_k,O_k;v_k)\|_{\mathcal{E}^{n/p}_{T^\ast_k}}\leq C\{C_0\eta+(M\eta)^{2}(1+M\eta)^{n+3}\}.
\end{equation}
By choosing $M=3CC_0$ and $\eta$ sufficient small enough such that
\be
 C(M\eta)(1+M\eta)^{n+3}\leq \f{1}{6},\nn
\ee
so \eqref{4.15} is followed by \eqref{4.155} directly.

Therefore, we obtain a sequence of approximate solutions $(\rho_k, F_k; u_k)$ to the system \eqref{1.1} on $[0,T_k)$ satisfying
\be\label{4.16}
\|(a_k,O_k;v_k)\|_{\mathcal{E}^{n/p}_{T_k}}\leq M\eta,
\ee
for any $k\in \bN$. From \eqref{4.8r} and \eqref{4.16}, we have
\bea
&&\|(a_k,O_k;v_k)\|_{E^{n/2}_{T_k}}\leq C\Big\{\|(a_{0,k},O_{0,k};v_{0,k})\|_{E^{n/2}_0}\nn\\
&&\qq\qq\qq\qq+\|(a_k,O_k;v_k)\|_{E^{n/2}_{T_k}}(M\eta)\big(1+M\eta\big)^{n+3}\Big\},\nn
\eea
which implies
\be\label{4.17}
\|(a_k,O_k;v_k)\|_{E^{n/2}_{T_k}}\leq C\|(a_{0,k},O_{0,k};v_{0,k})\|_{E^{n/2}_0},
\ee
where we chose $\eta$ sufficiently small. Consequently, based on Proposition \ref{pro4.2}, the continuity argument ensures that $T_k=+\i$ for any $k\in\bN$.
\subsection{Passing to the limit and existence}
Next, the existence of the solution will be proved by the compact argument.  We show that, up to an extraction, the sequence $(a_k,O_k;v_k)$ converges in the distributional sense to some function $(a,O;v)$ such that
\be\label{4.18}
(a,O;v)\in \Big(\t{L}^\i\B^{n/2-1,n/p}_{2,p}\cap{L}^1\B^{n/2+1,n/p}_{2,p}\Big)^{1+n^2}\times\Big(\t{L}^\i\B^{n/2-1,n/p-1}_{2,p}\cap{L}^1\B^{n/2+1,n/p+1}_{2,p}\Big).
^{n}
\ee
Indeed, it follows from \eqref{4.16} that $(a_k,O_k)$ is uniformly bounded in $\t{L}^\i(0,\i;\dot{B}^{n/p}_{p,1})$ and $v_k$ is uniformly bounded in $\t{L}^\i(0,\i;\dot{B}^{n/p-1}_{p,1})\cap L^1(0,\i;\dot{B}^{n/p+1}_{p,1})$. By interpolation, we also deduce that
$v_k$ is uniformly bounded in $\t{L}^{\f{2}{2-\ve}}(0,\i;\dot{B}^{n/p+1-\ve}_{p,1})$ for any $\ve\in [0,2]$. We claim that $(a_k,O_k; v_k)$ is uniformly bounded in
 \be\label{4.19x}
\Big(C^{1/2}_{loc}(\bR_+;\dot{B}^{n/p-1}_{p,1})\Big)^{1+n^2}\times\Big(C^{\f{2-\zeta}{2}}_{loc}(\bR_+;\dot{B}^{n/p-1-\zeta}_{p,1})\Big)^n
\ee
with $\zeta=\min\{\f{2n}{p}-1,1\}$, which is a direct consequence of
\be\label{4.19}
(\p_t a_k,\p_t O_k;\p_t v_k)\in \Big(\t{L}^2_{loc}\dot{B}^{n/p-1}_{p,1}\Big)^{1+n^2}\times \Big(\t{L}^{\f{2}{2-\zeta}}_{loc}\dot{B}^{n/p-1-\zeta}_{p,1}\Big)^n.
\ee
Recalling \eqref{r1}, we have
\be
\p_t a_k=-v_k\cdot\nabla a_k-\nabla\cdot v_k-a_k\nabla\cdot v_k \nn
\ee and
\be
\p_tO_k=-v_k\cdot\nabla O_k+\nabla v_k+\nabla v_kO_k.\nn
\ee

By interpolation and Lemma \ref{ll2.3}, it follows from \eqref{4.16} that
$(\p_t a_k,\p_t O_k) \in \Big(\t{L}^2_{loc}\dot{B}^{n/p-1}_{p,1}\Big)^{1+n^2}$,  which implies that $( a_k, O_k)$ is uniformly bounded in $\Big(C^{1/2}_{loc}(\bR_+;\dot{B}^{n/p-1}_{p,1})\Big)^{1+n^2}$. On the other hand,
\bea
&&\p_tv_k=-v_k\cdot\nabla v_k+\mathcal{A}v_k-\na a_k+\na\cdot O_k+O^{jl}_k\p_jO^{\bullet l}_k\nn\\
&&\qq -I(a_k)\mathcal{A}v_k-K(a_k)\na a_k+\f{1}{1+a_k}\text{div}\big(2\t{\mu}(a_k)D(v_k)+\t{\la}(a_k)\text{div} v_k\text{Id}\big).\nn
\eea
It's easy to see that
\be\label{4.20}
\|\mathcal{A}v_k\|_{\t{L}^{\f{2}{2-\zeta}}\dot{B}^{n/p-1-\zeta}_{p,1}}\ls \| v_k\|_{\t{L}^{\f{2}{2-\zeta}}\dot{B}^{n/p+1-\zeta}_{p,1}}.
\ee
Thanks to Lemma \ref{ll2.3} and Proposition \ref{pa.2}, we have
\bea\label{4.21}
&&\|\big(v_k\cdot\nabla v_k,I(a_k)\mathcal{A}v_k\big)\|_{\t{L}^{\f{2}{2-\zeta}}\dot{B}^{n/p-1-\zeta}_{p,1}}\nn\\
&\ls&\|v_k\|_{\t{L}^{2}\dot{B}^{n/p}_{p,1}}\|\nabla v_k\|_{\t{L}^{\f{2}{1-\zeta}}\dot{B}^{n/p-1-\zeta}_{p,1}}+\|a_k\|_{\t{L}^\i\dot{B}^{n/p}_{p,1}}\|\nabla^2 v_k\|_{\t{L}^{\f{2}{2-\zeta}}\dot{B}^{n/p-1-\zeta}_{p,1}}.
\eea
Also, due to the embedding $\B^{n/2-1,n/p}_{2,p}\hookrightarrow \dot{B}^{n/p-\zeta}_{p,1}$ and Proposition \ref{pa.2}, we arrive at
\bea\label{4.22}
&&\|(\na a_k,\na O_k)\|_{\t{L}^\i\dot{B}^{n/p-1-\zeta}_{p,1}}+\|(K(a_k)\na a_k,O_k\na O_k)\|_{\t{L}^\i\dot{B}^{n/p-1-\zeta}_{p,1}}\nn\\
&\ls&\big(1+\|(a_k,O_k)\|_{\t{L}^\i\dot{B}^{n/p}_{p,1}}\big)\|( a_k, O_k)\|_{\t{L}^\i\dot{B}^{n/p-\zeta}_{p,1}}.
\eea
As before, we write
 \be
 \f{1}{1+a_k}\text{div}\big(2\t{\mu}(a_k)D(v_k)\big)=\f{1}{1+a_k}\t{\mu}(a_k)\na^2 v_k+\f{1}{1+a_k}\na\t{\mu}(a_k)\na v_k.\nn
\ee
Then by applying Lemma \ref{ll2.3} and Proposition \ref{pa.2}, we get
\bea\label{e4.63ee}
&&\|\f{1}{1+a_k}\text{div}\big(2\t{\mu}(a_k)D(v_k)\big)\|_{\t{L}^{\f{2}{2-\zeta}}\dot{B}^{n/p-1-\zeta}_{p,1}}\nn\\
&\ls&\|\f{1}{1+a_k}\t{\mu}(a_k)\na^2 v_k\|_{\t{L}^{\f{2}{2-\zeta}}\dot{B}^{n/p-1-\zeta}_{p,1}}+\|\f{1}{1+a_k}\na\t{\mu}(a_k)\na v_k\|_{\t{L}^{\f{2}{2-\zeta}}\dot{B}^{n/p-1-\zeta}_{p,1}}\nn\\
&\ls&(1+\|I(a_k)\|_{\t{L}^\i\dot{B}^{n/p}_{p,1}})\Big(\|\t{\mu}(a_k)\na^2 v_k\|_{\t{L}^{\f{2}{2-\zeta}}\dot{B}^{n/p-1-\zeta}_{p,1}}+\|\na\t{\mu}(a_k)\na v_k\|_{\t{L}^{\f{2}{2-\zeta}}\dot{B}^{n/p-1-\zeta}_{p,1}}\Big)\nn\\
&\ls&(1+\|I(a_k)\|_{\t{L}^\i\dot{B}^{n/p}_{p,1}})\Big(\|\t{\mu}(a_k)\|_{\t{L}^\i\dot{B}^{n/p}_{p,1}}\|\na^2 v_k\|_{\t{L}^{\f{2}{2-\zeta}}\dot{B}^{n/p-1-\zeta}_{p,1}}\nn\\
&&\qq\qq\qq\qq\qq\qq +\|\na\t{\mu}(a_k)\|_{\t{L}^\i\dot{B}^{n/p-1}_{p,1}}\|\na v_k\|_{\t{L}^{\f{2}{2-\zeta}}\dot{B}^{n/p-\zeta}_{p,1}}\Big)
\eea
and $\f{1}{1+a_k}\text{div}\big(\t{\la}(a_k)\text{div} v_k\text{Id}\big)$ may be treated along the same way.  Consequently,
combining $\eqref{4.20}-\eqref{e4.63ee}$, we conclude that
\be
\p_t v_k\in  \Big(L^{\f{2}{2-\zeta}}_{loc}\dot{B}^{n/p-1-\zeta}_{p,1}\Big)^n,\nn
\ee
which implies that $v_k$ is uniformly bounded in $C^{\f{2-\zeta}{2}}_{loc}(\bR_+;\dot{B}^{n/p-1-\zeta}_{p,1})$. Therefore the claim \eqref{4.19x} is proved.
Furthermore, we see that $(a_k,O_k;v_k)$ is equicontinuous on $\bR_+$ valued in $\big(\dot{B}^{n/p-1}_{p,1}\big)^{1+n^2}\times \big(\dot{B}^{n/p-1-\zeta}_{p,1}\big)^n$. Let $\{\phi_j\}_{j\in\bN}$ be a sequence of smooth functions supported in  the ball $B(0,j+1)$ and equal to $1$ in $B(0,j)$. It follows from  \eqref{4.19x} that $(\phi_ja_k,\phi_jO_k;\phi_jv_k)$ is uniformly bounded in
\be
\Big(C^{1/2}_{loc}(\bR_+;\dot{B}^{n/p-1}_{p,1})\Big)^{1+n^2}\times\Big(C^{\f{2-\zeta}{2}}_{loc}(\bR_+;\dot{B}^{n/p-1-\zeta}_{p,1})\Big)^n.\nn
\ee
Observe that the map $(a_k,O_k;v_k)\mapsto (\phi_ja_k,\phi_jO_k;\phi_jv_k)$ is compact from
\be
\big(\dot{B}^{n/p-1}_{p,1}\cap\dot{B}^{n/p}_{p,1}\big)^{1+n^2}\times\big(\dot{B}^{n/p-1-\zeta}_{p,1}\cap\dot{B}^{n/p-1}_{p,1}\big)^n\nn
\ee
into
\be
\big(\dot{B}^{n/p-1}_{p,1}\big)^{1+n^2}\times\big(\dot{B}^{n/p-1-\zeta}_{p,1}\big)^n.\nn
\ee
By applying Ascoli's theorem and Cantor's diagonal process, there exist a $(a,O;v)$ such that for any smooth function $\phi\in C^\i_0(\bR^n)$,
\bea\label{4.24}
&&(\phi a_k,\phi O_k)\rightarrow (\phi a,\phi O) \q \text{in}\q \big(L^\i(\bR_+;\dot{B}^{n/p-1}_{p,1})\big)^{1+n^2},\nn\\
&&\phi v_k\rightarrow \phi v \q \text{in}\q \big(L^\i(\bR_+;\dot{B}^{n/p-1-\zeta}_{p,1})\big)^n,
\eea
when $k\rightarrow +\i$ (up to an extraction). Actually, by interpolation, we also have
 \bea\label{4.25}
&&(\phi a_k,\phi O_k)\rightarrow (\phi a,\phi O) \q \text{in}\q \big(L^\i(\bR_+;\dot{B}^{n/p-s}_{p,1})\big)^{1+n^2}\q \forall\ 0<s\leq 1,\nn\\
&&\phi v_k\rightarrow \phi v \q \text{in}\q \big(L^1(\bR_+;\dot{B}^{n/p+s}_{p,1})\big)^n\q \forall -1\leq s<1.
\eea
Then, using the so-called Fatou property in Besov spaces and the uniform bound in \eqref{4.16}, we conclude that \eqref{4.18} is fulfilled.
It is a routine process to verify that $(a,O; v)$ satisfies the system \eqref{r1} in the sense of distributions. Below is to check the desired regularity of solutions. Noticing that
\bea
&&\p_t a+v\cdot\nabla a=-\nabla\cdot v-a\nabla\cdot v\in {L}^1_{loc}(\B^{n/2-1,n/p}_{2,p})\cap L^1(\B^{n/2+1,n/p}_{2,p}), \nn\\
&&\p_tO+v\cdot\nabla O=\nabla v+\nabla vO\in {L}^1_{loc}(\B^{n/2-1,n/p}_{2,p})\cap L^1(\B^{n/2+1,n/p}_{2,p}).\nn
\eea
since $(a_0,O_0)\in\B^{n/2-1,n/p}_{2,p}$, the classical result for transport equations indicates that
\be
 (a,O)\in C(\bR_+;\B^{n/2-1,n/p}_{2,p}).\nn
 \ee
On the other hand,
\bea
&&\p_tv-\mathcal{A}v=-v\cdot\nabla v-\na a+\na\cdot O+O^{jk}\p_jO^{\bullet k}-I(a)\mathcal{A}v-K(a)\na a,\nn\\
&&\qq\qq\qq+\f{1}{1+a}\text{div}\big(2\t{\mu}(a)D(v)+\t{\la}(a)\text{div} v\text{Id}\big),\nn\\
&&\qq\qq\in \t{L}^1_{loc}(\B^{n/2-1,n/p-1}_{2,p}),\nn
\eea
So the maximal regularity of heat equation enables us to get $v\in C(\bR_+; \B^{n/2-1,n/p-1}_{2,p}).$
\subsection{Uniqueness}
Due to techniques, allow us to only deal with the case $2\leq p\leq n$ in the proof uniqueness of solutions. We will consider the remaining interval with respect to $p$ in near future. The proof depends on a logarithmic inequality. For convenience of reader, we present it by a lemma.
\begin{lemma}[\cite{Dr4}]\label{l4.2}
Let $s\in \bR$. Then for any $1\leq p,r\leq +\i$ and $0<\ve\leq 1$, we have
\be
\|f\|_{\t{L}^r_t\dot{B}^s_{p,1}}\leq C\f{\|f\|_{\t{L}^r_t\dot{B}^s_{p,\i}}}{\ve}\log \Big(e+\f{\|f\|_{\t{L}^r_t\dot{B}^{s-\ve}_{p,\i}}+\|f\|_{\t{L}^r_t\dot{B}^{s+\ve}_{p,\i}}}{\|f\|_{\t{L}^r_t\dot{B}^s_{p,\i}}}\Big).\nn
\ee
\end{lemma}

Assume that $(\rho_i,F_i;u_i)(i=1,2)$ are two solution to the system \eqref{1.1} with the same initial data. Without loss of generality, we may assume that
\be\label{4.26}
\|(\rho_i-1,F_i-I;u_i)\|_{\mathcal{E}^{n/p}}\leq M\eta. \q \text{for}\q i=1,2.
\ee
Using embedding and \eqref{4.26}, we have
\be
\|\rho_i-1\|_{L^\i(\bR_+\times\bR^n)}\leq C\|\rho_i-1\|_{\mathcal{E}^{n/p}}\leq CM\eta\leq \f{1}{2},\q \text{for}\q i=1.2\nn
\ee
for $\eta>0$ sufficiently small. Set
\bea
&&a_i(t,x)=\rho_i(\chi^2_0t,\chi_0x)-1, \nn\\
&&O_i(t,x)=F_i(\chi^2_0t,\chi_0x)-I,\nn \\
&&v_i(t,x)=\chi_0u_i(\chi^2_0t,\chi_0x),\nn
\eea
for $i=1,2$.
and
\be
\dl a =a_1-a_2, \dl O=O_1-O_2; \dl v=v_1-v_2.\nn
\ee
Thanks to \eqref{r1}, we find that $(\dl a,\dl v,\dl O)$ satisfies
\be\label{4.27}
\left\{
\begin{aligned}
&\p_t \dl a+v_2\cdot\nabla \dl a=\dl F, \\
&\p_t\dl v-\mathcal{A}\dl v=\dl G,\\
&\p_t\dl O+v_2\cdot\nabla \dl O=\dl H,\\
&(\dl a,\dl O;\dl v)=(0,0,0),
\end{aligned}
\right.
\ee
with
\bea
&&\dl F=-\dl v\cdot\na a_1-\na\cdot\dl v-a_1\na\cdot\dl v-\dl a\na\cdot v_2,\nn\\
&&\dl H=\dl v\cdot\na O_1+\na \dl v+\na\dl vO_1+\na v_2\dl O,\nn\\
&&\dl G=-\na \dl a+\na\cdot \dl O-(v_1\cdot\na v_1-v_2\cdot\na v_2)+(O^{jk}_1\p_jO^{\bullet k}_1-O^{jk}_2\p_jO^{\bullet k}_2)\nn\\
&&\qq\q -I(a_1)\mathcal{A}v_1+I(a_2)\mathcal{A}v_2-K(a_1)\na a_1+K(a_2)\na a_2\nn\\
&&\qq\q+\f{1}{1+a_1}\text{div}\big(2\t{\mu}(a_1)D(v_1)+\t{\la}(a_1)\text{div} v_1\text{Id}\big)\nn\\
&&\qq\q-\f{1}{1+a_2}\text{div}\big(2\t{\mu}(a_2)D(v_2)+\t{\la}(a_2)\text{div} v_2\text{Id}\big).
\eea
In the following, we denote
\be
V_i(t)=\int^t_0\|v_i(\tau)\|_{\dot{B}^{n/p+1}_{p,1}}d\tau\q \text{for}\ i=1,2
\ee
and we denote by $A_t$ a constant depending on $\|a_i\|_{\t{L}^\i_t\dot{B}^{n/p}_{p,1}}$ for $i=1,2$. Due to the embedding $\mathcal{E}^{n/p}\subseteq\mathcal{E}^1(p\leq n)$, it is suffices to prove the uniqueness in $\mathcal{E}^1$. So we
take $p=n$ in the subsequent process.

Applying Proposition \ref{pa.3}, we get
\be\label{4.30}
\|\big(\dl a(t),\dl O(t)\big)\|_{\dot{B}^0_{p,\i}}\leq e^{CV_2(t)}\int^t_0\|\big(\dl F(\tau),\dl H(\tau)\big)\|_{\dot{B}^0_{p,\i}}d\tau.
\ee
By Lemma \ref{ll2.3}, we have
\bea
&&\|\big(\dl F(\tau),\dl H(\tau)\big)\|_{\dot{B}^0_{p,\i}}\nn \\
&\ls& \|v_2\|_{\dot{B}^2_{p,1}}\|(\dl a,\dl O) \|_{\dot{B}^0_{p,\i}}+\big(1+\|(a_1,O_1)\|_{\dot{B}^1_{p,1}}\big)\|\dl v\|_{\dot{B}^1_{p,1}}.\nn
\eea
Inserting the equality into \eqref{4.30}, we arrive at by Gronwall's inequality
\be\label{4.31}
\|\big(\dl a(t),\dl O(t)\big)\|_{\dot{B}^0_{p,\i}}\leq e^{CV_2(t)}\int^t_0\big(1+\|(a_1,O_1)\|_{\dot{B}^1_{p,1}}\big)\|\dl v\|_{\dot{B}^1_{p,1}}d\tau.
\ee
Using Proposition \ref{pA.4} to the second equation of \eqref{4.27} gives
\be\label{4.32}
\|\dl v\|_{\t{L}^1_t\dot{B}^1_{p,\i}}+\|\dl v\|_{\t{L}^2_t\dot{B}^0_{p,\i}}\ls \|\dl G(\tau)\|_{\t{L}^1_t\dot{B}^{-1}_{p,\i}}.
\ee
Furthermore, by Lemma \ref{ll2.3} and Proposition \ref{pa.2}, it is shown that
\bea\label{4.33}
&&\|\dl G(\tau)\|_{\t{L}^1_t\dot{B}^{-1}_{p,\i}}\ls \|(v_1,v_2)\|_{\t{L}^2_t\dot{B}^{1}_{p,1}}\|\dl v\|_{\t{L}^2_t\dot{B}^{0}_{p,\i}}+A_t\|a_1\|_{\t{L}^\i_t\dot{B}^{1}_{p,1}}\|\dl
v\|_{\t{L}^1_t\dot{B}^{1}_{p,\i}} \nn\\
&&\qq\qq\qq +A_t\int^t_0(1+\|v_2\|_{\dot{B}^{2}_{p,1}})\|(\dl a,\dl O)\|_{\dot{B}^{0}_{p,\i}}d\tau.
\eea
According to our a prior estimates, by choosing $\eta$ small, we have
\be
A_t\|a_1\|_{\t{L}^\i_t\dot{B}^{1}_{p,1}}+\|(v_1,v_2)\|_{\t{L}^2_t\dot{B}^{1}_{p,1}}\ls M\eta\ll 1.\nn
\ee
Consequently, inserting \eqref{4.33} into \eqref{4.32} to implies that
\bea\label{4.34}
\|\dl v\|_{\t{L}^1_t\dot{B}^1_{p,\i}}+\|\dl v\|_{\t{L}^2_t\dot{B}^0_{p,\i}}&\ls& A_t\int^t_0(1+\|v_2\|_{\dot{B}^{2}_{p,1}})\|(\dl a,\dl O)\|_{\dot{B}^{0}_{p,\i}}d\tau.
\eea
Combining \eqref{4.31} and \eqref{4.34}, we get
\be\label{4.74ee}
\|\dl v\|_{\t{L}^1_t\dot{B}^1_{p,\i}}\ls \int^t_0(1+\|v_2\|_{\dot{B}^{2}_{p,1}})\|\dl v\|_{\t{L}^1_\tau\dot{B}^1_{p,1}}d\tau.
\ee
Applying Lemma \ref{l4.2} with $s=r=\ve=1$ and $f=\dl v$, we obtain
\be
\|\dl v\|_{\t{L}^1_t\dot{B}^1_{p,1}}\leq C\|\dl v\|_{\t{L}^1_t\dot{B}^1_{p,\i}}\log \Big(e+\f{\|\dl v\|_{\t{L}^1_t\dot{B}^{0}_{p,\i}}+\|\dl v\|_{\t{L}^1_t\dot{B}^{2}_{p,\i}}}{\|\dl v\|_{\t{L}^1_t\dot{B}^1_{p,\i}}}\Big),\nn
\ee
which together with \eqref{4.34} and \eqref{4.31} indicates that
\bea
&&\|\dl v(t)\|_{\t{L}^1_t\dot{B}^1_{p,\i}}\leq e^{CV_2(t)}A_t\int^t_0(1+\|v_2\|_{\dot{B}^{2}_{p,1}})\|\dl v\|_{\t{L}^1_\tau\dot{B}^{1}_{p,\i}}\nn\\
&&\qq\qq \qq\qq \times\log\Big(e+C_\tau\|\dl v\|^{-1}_{\t{L}^1_\tau\dot{B}^{1}_{p,\i}}\Big) d\tau,\nn
\eea
where $C_\tau=\|\dl v\|_{\t{L}^1_\tau\dot{B}^{0}_{p,\i}}+\|\dl v\|_{\t{L}^1_\tau\dot{B}^{2}_{p,\i}}$. Noting $\|v_2\|_{\dot{B}^{2}_{p,1}}$ is integrable on $[0,\i]$ and
\be
\int^1_0\f{dr}{r\log(r+C_tr^{-1})}=+\i,
\ee
the Osgood lemma implies that $(\dl a,\dl O;\dl v)=0$ on $[0,t]$. So a continuity argument ensures that $(a_1,O_1;v_1)=(a_2,O_2;v_2)$  for any $t\in[0,\i)$.

\section{The proof of time-decay estimates}
In this section, we aim at proving the time-weighted energy inequality \eqref{de1} taking for granted Theorem \ref{thglobal2}. We will
proceed the proof into the three subsections, according to three terms of the definition of $\mathcal{G}_p(t)$.
Subsection 5.1 is devoted to the low-frequency estimates. In the spirit of \cite{DX1}, we need to perform nonlinear estimates in terms of deformation tensor in the Besov space with negative regularity. In Subsection 5.2, in order to overcome the technical difficulty that there is loss of one derivative for the density and deformation tensor at high frequencies, we develop ``\emph{two effective velocities}'' to obtain the upper bound for the second term in $\mathcal{G}_p(t)$.
To close the high-frequency estimates, in Subsection 5.3, a crucial observation enables us to establish gain of regularity and decay altogether for the velocity, which strongly depends on Proposition \ref{pA.4}.

For simplicity, we denote
\be\label{e6.1}
\mathcal{X}_p(t)\triangleq \|(a,O;v)\|_{\E^{n/p}_t}.
\ee
In what follows, we will use the two key lemmas repeatedly.
\begin{lemma}\label{l6.1}
Let $0\leq \sigma_1\leq\sigma_2$ with $\sigma_2>1$. It holds that
\be
\int_{0}^{t}\langle t-\tau\rangle^{-\sigma_{1}}\langle\tau\rangle^{-\sigma_{2}}d\tau\lesssim\langle t\rangle^{-\sigma_{1}}
\ee
and
\be
\int^t_0 \langle t-\tau\rangle^{-\sigma_1}\tau^{-\th}\langle \tau\rangle^{\th-\sigma_2}d\tau\ls \langle t\rangle^{-\sigma_1}\q \text{if}\ 0\leq \th<1.
\ee
\end{lemma}

\begin{lemma}[\cite{DX1}]\label{lxx}
Let $X:[0,T]\rightarrow \bR_+$ be a continuous function such that $X^p$ is a differentiable for some $p\geq 1$ and satisfies
\be
\f{1}{p}\f{d}{dt}X^p+BX^p\leq AX^{p-1}\nn
\ee
for some constant $B\geq 0$ and measurable function $A: [0,T]\rightarrow \bR_+$. Define $X_\dl=(X^p+\dl^p)^{\f{1}{p}}$ for $\dl>0$. Then it holds that
\be\label{6.2e}
\f{d}{dt}X_\dl+BX_\dl\leq A+B\dl.
\ee
\end{lemma}
For convenience, we denote by $\|\cdot\|_{\dl,L^p}:=(\|\cdot\|^p_{L^p}+\dl^p)^{1/p}$ for $1\leq p<\i$.

\subsection{Bounds for the low frequencies}
From \eqref{e4.10} and \eqref{e4.11}, we have
\bea
\f{d}{dt}\big(\|(a_k,O_k;v_k)\|^2_{L^2}\big)+2^{2k}\|(a_k,O_k;v_k)\|^2_{L^2}
\lesssim \big(\sum\limits_{i=0,1,3,4}\|\dD_k G_i\|_{L^2}\big)\|(a_k,O_k;v_k)\|_{L^2}.\nn
\eea
It follows from Lemma \ref{lxx} that
\bea
\f{d}{dt}\|(a_k,O_k;v_k)\|_{\dl,L^2}+2^{2k}\|(a_k,O_k;v_k)\|_{\dl,L^2}\lesssim \sum\limits_{i=0,1,3,4}\|\dD_k G_i\|_{L^2}+2^{2k}\dl.\nn
\eea
Then integrating in time and letting $\dl\rightarrow 0$, then there exists a $c_0$ such that
\bea\label{6.1}
\|(a_k,O_k;v_k)\|_{L^2}&\ls& e^{-c_02^{2k}t}\|(\dD_ka_0,\dD_kO_{0};\dD_kv_{0})\|_{L^2}\nn\\
&&+\int^t_0e^{c_02^{2k}(\tau-t)}\sum\limits_{i=0,1,3,4}\|\dD_k G_i\|_{L^2}d\tau.
\eea
Regarding the first term in \eqref{6.1}, we multiply the factor $\langle t\rangle ^{\f{s+s_0}{2}}2^{ks}$ and sum up on $2^k\leq R_0$ to get
\bea\label{e6.2}
&&\langle t\rangle^{\f{s+s_0}{2}}\sum\limits_{2^k\leq R_0}2^{ks}e^{-c_02^{2k}t}\|(\dD_ka_{0},\dD_kO_{0};\dD_kv_{0})\|_{L^2}\nn\\
&\ls&\|(a_{0},O_{0};v_{0})\|^\ell_{\dot{B}^{-s_0}_{2,\i}}\sum\limits_{2^k\leq R_0}(2^k\langle t\rangle^{\f{1}{2}})^{s+s_0}e^{-c_0(2^k\s{t})^2}.\nn\\
&\ls&\|(a_{0},O_{0};v_{0})\|^\ell_{\dot{B}^{-s_0}_{2,\i}}\Big(\sum\limits_{2^k\leq R_0}(2^k\s{t})^{s+s_0}e^{-c_0(2^k\s{t})^2}+2^{k_0(s+s_0)}\Big)\nn\\
&\ls& \|(a_{0},O_{0};v_{0})\|^\ell_{\dot{B}^{-s_0}_{2,\i}}, \label{3.4}
\eea
where we have used the fact $\sum\limits_{2^k\leq R_0}(2^k\s{t})^{s+s_0}e^{-c_0(2^k\s{t})^2}\leq C$ when $s+s_0> 0$.
So we have
\bea\label{6.3}
\sum\limits_{2^k\leq R_0}2^{ks}e^{-c_02^{2k}t}\|(\dD_ka_{0},\dD_kO_{0};\dD_kv_{0})\|_{L^2}&\ls& \langle t\rangle^{-\f{s+s_0}{2}}\|(a_{0},O_{0};v_{0})\|^\ell_{\dot{B}^{-s_0}_{2,\i}}.
\eea
Furthermore, the corresponding nonlinear term in \eqref{6.1} can be estimated as
\bea
\sum\limits_{2^k\leq R_0}2^{ks}\int^t_0e^{c_02^{2k}(\tau-t)}\sum\limits_{i=0,1,3,4}\|\dD_k G_i\|_{L^2}d\tau
\lesssim \int^t_0 \langle t-\tau\rangle^{-\f{s+s_0}{2}}\sum\limits_{i=0,1,3,4}\| G_i\|^\ell_{\dB^{-s_0}_{2,\i}}d\tau.\label{3.8}
\eea
We claim that if $p$ fulfills the assumption as in Theorem \ref{thdecay}, then we have for all $t\geq0$,
\bea
\int^t_0 \langle t-\tau\rangle^{-\f{s+s_0}{2}}\sum\limits_{i=0,1,2,3}\| G_i\|^\ell_{\dB^{-s_0}_{2,\i}}d\tau
\lesssim \langle t\rangle^{-\f{s+s_0}{2}}\big(\mathcal{G}^2_p(t)+\mathcal{X}^2_p(t)\big),   \label{3.9}
\eea
where $\mathcal{G}_p(t)$ and $\mathcal{X}_p(t)$ are defined by \eqref{dfD} and \eqref{e6.1}.

Since those quadratic terms containing $a$ and $v$ in $G_i(i=0,1,3,4)$ have already been done in \cite{DX1}, it suffices to give suitable decay estimates for some terms involving in $O$. Precisely, we need to hand the following integral
\bea
&&\int^t_0 \langle t-\tau\rangle^{-\f{s+s_0}{2}}\|\big(\p_i(aO^{ij}),\p_kv^iO^{kj},v\cdot\na O^{ij},\nn\\
&&\qq\qq\qq O^{jk}\p_jO^{ik} ,O^{lj}\p_lO^{ik},O^{lk}\p_lO^{ij}\big)\|^\ell_{\dB^{-s_0}_{2,\i}}d\tau.\nn
\eea
As far as we know, the regularity level remains the same between the density and deformation tensor. Hence, these terms $\p_iaO^{ij}, a\p_iO^{ij}$,
$O^{jk}\p_jO^{ik}$, $O^{lj}\p_lO^{ik}$, $O^{lk}\p_lO^{ij}$ can be treated along the same line. In principle, the above integral can be reduced to
\bea\label{6.7}
&&\int^t_0 \langle t-\tau\rangle^{-\f{s+s_0}{2}}\|\big(\na a\cdot O,\na vO,v\cdot\na O \big)\|^\ell_{\dB^{-s_0}_{2,\i}}d\tau.
\eea
We decompose \eqref{6.7} as follows
\be
\eqref{6.7}\triangleq I^\ell+I^h,\nn
\ee
where
\be
I^\ell=\int^t_0 \langle t-\tau\rangle^{-\f{s+s_0}{2}}\|\big(O\cdot\na a^\ell,O\na v^\ell,v \cdot\na O^\ell\big)\|^\ell_{\dB^{-s_0}_{2,\i}}d\tau,\nn
\ee
and
\be
I^h=\int^t_0 \langle t-\tau\rangle^{-\f{s+s_0}{2}}\|\big(O\cdot\na a^h,O\na v^h,v\cdot\na O^h\big)\|^\ell_{\dB^{-s_0}_{2,\i}}d\tau.\nn
\ee

In order to handle $I^\ell$ in terms of with $a^\ell, O^\ell$ and $v^\ell$, we use the following Lemma.
\begin{lemma}\label{lemmaaa}
Let $s_0=n(2/p-1/2)$ and $p$ satisfy the assumption in  Theorem \ref{thdecay}. It holds that
\be
\|fg\|_{\dot{B}^{-s_0}_{2,\i}}\ls\|f\|_{\dot{B}^{1-n/p}_{p,1}}\|g\|_{\dot{B}^{n/2-1}_{2,1}},  \label{3.10}
\ee
and
\be
\|fg\|_{\dot{B}^{-n/p}_{2,\i}}\ls\|f\|_{\dot{B}^{n/p-1}_{p,1}}\|g\|_{\dot{B}^{1-n/p}_{2,1}}.  \label{3.11}
\ee
\end{lemma}
The reader is referred to \cite{DX1} for the detailed proof. Owing to the embedding theorem and the definition of $\mathcal{G}_p(t)$, we shall often use the following inequalities
\be
\|(a,O;v)^\ell(\tau)\|_{\dot{B}^{1-\f{n}{p}}_{p,1}}\ls \|(a,O;v)^\ell(\tau)\|_{\dot{B}^{1-s_0}_{2,1}}\ls\langle \tau\rangle^{-\f{1}{2}}\mathcal{G}_p(\tau) ,   \label{3.12}
\ee
and
\be
\|(a,O)\|_{\dot{B}^{\f{n}{p}}_{p,1}}\ls \langle \tau\rangle^{-\f{n}{p}}\mathcal{G}_p(\tau).  \label{3.13}
\ee
Indeed, the above inequality is obvious for the high frequencies since $\alpha \geq \frac{n}{p}$, and we have
\bea\label{6.12}
\|(a,O)^\ell\|_{\dot{B}^{\f{n}{p}}_{p,1}}&\ls&\|(a,O)\|^\ell_{\dot{B}^{\f{n}{2}}_{2,1}}\ls\langle \tau\rangle^{-\f{1}{2}(s_0+n/2)}\mathcal{G}_p(\tau)=\langle \tau\rangle^{-\f{n}{p}}\mathcal{G}_p(\tau).
\eea
Notice that  $1-\f{n}{p}\leq \f{n}{p}$ and the definition of $\mathcal{G}_p(t)$, we arrive at
\bea\label{3.14}
\|v^h\|_{\dot{B}^{1-\f{n}{p}}_{p,1}}\ls\|v^h\|_{\dot{B}^{\f{n}{p}}_{p,1}}\ls\bigg(\|v^h\|_{\dot{B}^{\f{n}{p}-1}_{p,1}}\|\nabla v^h\|_{\dot{B}^{\f{n}{p}}_{p,1}}\bigg)^{\f{1}{2}}
            \ls\tau^{-\f{1}{2}}\langle \tau\rangle^{-\f{\al}{2}}\mathcal{G}_p(\tau).
\eea
Next, we begin with bound $I^\ell$ and $I^h$.

\noindent\textbf{Estimates for $I^\ell$}\\
Taking advantage of \eqref{3.10}, \eqref{3.12}, \eqref{3.13} and \eqref{3.14}, we get
\bea
\int^t_0\langle t-\tau\rangle^{-\f{s+s_0}{2}}\|(v\cdot\nabla O^\ell)\|^\ell_{\dot{B}^{-s_0}_{2,\i}}d\tau&\ls&\int^t_0\langle t-\tau\rangle^{-\f{s+s_0}{2}}\|v\|_{\dot{B}^{1-\f{n}{p}}_{p,1}}\|\nabla O^\ell\|_{\dot{B}^{\f{n}{2}-1}_{2,1}}d\tau\nn\\
                &\ls& \mathcal{G}^2_p(t)\int^t_0\langle t-\tau\rangle^{-\f{s+s_0}{2}}\Big(\langle \tau\rangle^{-\f{1}{2}}+\tau^{-\f{1}{2}}\langle \tau\rangle^{-\f{\alpha}{2}}\Big)\langle \tau\rangle^{-\f{n}{p}}d\tau. \nn
\eea
Due to the fact that $\f{n}{p}+\f{1}{2}>1$ and $\f{s+s_0}{2}\leq \f{n}{p}+\f{1}{2}$ for all $s\leq 1+\f{n}{2}$,  Lemma \eqref{l6.1} implies that
\be
\int^t_0\langle t-\tau\rangle^{-\f{s+s_0}{2}}\|(v\cdot\nabla O^\ell)\|^\ell_{\dot{B}^{-s_0}_{2,\i}}d\tau\ls \langle t\rangle^{-\f{s+s_0}{2}}\mathcal{G}^2_p(t).  \label{3.15}
\ee
The terms $O\cdot\nabla a^\ell$ and $O\na v^\ell$ can be treated along with the same lines, so we feel free to skip the details.

\noindent\textbf{Estimates for $I^h$}\\
For the term $I^h$ containing $a^h, O^h$ and $v^h$, as in \cite{DX1}, we proceed differently depending on whether $p>n$ and $p\leq n$. Let's first consider the easy case $2\leq p\leq n$. Applying \eqref{e2.2} with $\sigma=\f{n}{p}-1$ yields
\bea
\|fg^h\|_{\dot{B}^{-s_0}_{2,\i}}\ls\|f\|_{\dot{B}^{1-\f{n}{p}}_{p,1}}\Big(\|\dot{S}_{k_0+N_0}g^h\|_{L^{p^*}}+\|g^h\|_{\dot{B}^{\f{n}{p}-1}_{p,1}}\Big)
                                 \ls\|f\|_{\dot{B}^{1-\f{n}{p}}_{p,1}}\|g^h\|_{\dot{B}^{\f{n}{p}-1}_{p,1}}, \label{3.17}
\eea
where we have used the Berstein inequality ($p^\ast=\f{2p}{p-2}\geq p$) and the fact that only finite middle frequencies of $g$ are involving in $\dot{S}_{k_0+N_0}g^h$.\footnote{The limit case $p=n$ follows from $\|fg^h\|^\ell_{\dB^{-s_0}_{2,\i}}\ls \|fg^h\|_{L^{\f{n}{2}}}\ls \|f\|_{L^n}\|g^h\|_{L^n}\ls\|f\|_{\dB^0_{n,1}}\|g^h\|_{\dB^0_{n,1}}$.}

Taking $f=v$ and $g=\nabla O$ in \eqref{3.17}, we get
\bea
&&\int^t_0\langle t-\tau\rangle^{-\f{s+s_0}{2}}\|v\cdot\nabla O^h\|^\ell_{\dot{B}^{-s_0}_{2,\i}}d\tau\nn\\
&\ls& \int^t_0\langle t-\tau\rangle^{-\f{s+s_0}{2}}\|v\|_{\dot{B}^{1-\f{n}{p}}_{p,1}}\|\nabla O^h\|_{\dot{B}^{\f{n}{p}-1}_{p,1}}d\tau. \label{3.18}
\eea
It follows from $\eqref{3.12}$ and $\eqref{3.14}$ that
\be
\|v\|_{\dot{B}^{1-\f{n}{p}}_{p,1}}\ls \big(\langle \tau\rangle^{-\f{1}{2}}+\tau^{-\f{1}{2}}\langle\tau\rangle^{-\f{\al}{2}}\big)\mathcal{G}_p(\tau). \label{3.19}
\ee
The definition of $\mathcal{G}_p(t)$ implies that
\be
\|\nabla O^h\|_{\dot{B}^{\f{n}{p}-1}_{p,1}}\ls \langle\tau\rangle^{-\alpha}\mathcal{G}_p(\tau)\ \ \ \mbox{with}\ \ \ \alpha=\frac{n}{p}+\frac{1}{2}-\varepsilon. \label{3.20}
\ee
Inserting $\eqref{3.19}$ and \eqref{3.20} into \eqref{3.18}, we conclude that for $-s_0\leq s\leq \frac{d}{2}+1$,
\bea\label{3.21}
&&\int^t_0\langle t-\tau\rangle^{-\f{s+s_0}{2}}\|v\cdot\nabla O^h\|^\ell_{\dot{B}^{-s_0}_{2,\i}}d\tau\nn\\
&\ls&\mathcal{G}^2_p(t) \int^t_0\langle t-\tau\rangle^{-\f{s+s_0}{2}}\langle\tau\rangle^{-\alpha}\big(\langle \tau\rangle^{-\f{1}{2}}+\tau^{-\f{1}{2}}\langle\tau\rangle^{-\f{\al}{2}}\big)d\tau\nn\\
&\ls&\langle t\rangle^{-\f{s+s_0}{2}}\mathcal{G}^2_p(t).
\eea
Handling with the term $O\cdot\na a^h$ is similar. With aid of \eqref{3.17}, we have
\bea\label{3.22}
\int^t_0\langle t-\tau\rangle^{-\f{s+s_0}{2}}\|O\cdot\na a^h\|^\ell_{\dot{B}^{-s_0}_{2,\i}}d\tau
&\ls&\int^t_0\langle t-\tau\rangle^{-\f{s+s_0}{2}}\|O\|_{\dot{B}^{1-\f{n}{p}}_{p,1}}\|\nabla a^h\|_{\dot{B}^{\f{n}{p}-1}_{p,1}}d\tau\nn\\
&\ls&\mathcal{G}^2_p(t)\int^t_0\langle t-\tau\rangle^{-\f{s+s_0}{2}}\langle\tau\rangle^{-\f{1}{2}-\al}d\tau\nn\\
&\ls&\langle t\rangle^{-\f{s+s_0}{2}}\mathcal{G}^2_p(t).
\eea
Regarding the term $O\nabla v^h$, combining the embedding $L^\f{p}{2}\hookrightarrow \dot{B}^{-s_0}_{2,\i}$ and H\"{o}lder inequality, we obtain
\bea\label{3.23}
\int^t_0\langle t-\tau\rangle^{-\f{s+s_0}{2}}\|O\nabla v^h\|^\ell_{\dot{B}^{-s_0}_{2,\i}}d\tau
\ls\int^t_0\langle t-\tau\rangle^{-\f{s+s_0}{2}}\|O(\tau)\|_{L^p}\|\nabla v^h(\tau)\|_{L^p}d\tau.
\eea
By embedding, the definition of $\mathcal{G}_p(t)$ and the fact that $\alpha\geq\frac{n}{2p}$ for sufficiently small $\varepsilon>0$, we have
\bea\label{3.24}
\|O\|_{L^p}&\ls&\|O^\ell\|_{L^p}+\|O^h\|_{L^p}\ls\|O\|^\ell_{\dot{B}^{\f{n}{2}-\f{n}{p}}_{2,1}}+\|O\|^h_{\dot{B}^{\f{n}{p}-1}_{p,1}}\nn\\
             &\ls&(\langle \tau\rangle^{-\f{n}{2p}}+\langle \tau\rangle^{-\al})\mathcal{G}_p(t)\ls\langle \tau\rangle^{-\f{n}{2p}}\mathcal{G}_p(t).
\eea
Arguing as for proving \eqref{3.14}, it is easy to get for $2\leq p\leq n$,
\bea\label{3.25}
\|\nabla v^h(\tau)\|_{L^p}\ls\| v^h(\tau)\|_{\dot{B}^{\f{n}{p}}_{p,1}}\ls \tau^{-\f{1}{2}}\langle \tau\rangle^{-\f{\al}{2}}\mathcal{G}_p(\tau).
\eea
Furthermore, together with \eqref{3.24}-\eqref{3.25},  we have
\bea\label{3.26}
&&\int^t_0\langle t-\tau\rangle^{-\f{s+s_0}{2}}\|O\nabla v^h\|^\ell_{\dot{B}^{-s_0}_{2,\i}}d\tau\nn\\
&\ls&\mathcal{G}^2_p(t)\int^t_0\langle t-\tau\rangle^{-\f{s+s_0}{2}}\tau^{-\f{1}{2}}\langle \tau\rangle^{-(\f{\al}{2}+\f{n}{2p})}d\tau\ls\mathcal{G}^2_p(t)\langle t\rangle^{-\f{s+s_0}{2}}.
\eea
Let's end that step by considering $I^h$  involving $a^h, O^h$ and $v^h$ in the case of $p>n$.
Applying Inequality \eqref{e2.1} with $\sigma=1-\f{n}{p}$ and the embedding $\dot{B}^{\f{n}{p}}_{2,1}\hookrightarrow L^{p*}$ give that
\bea\label{3.28}
\|fg^h\|^\ell_{\dot{B}^{-s_0}_{2,\i}}&\ls& \big(\|f\|_{\dot{B}^{1-\f{n}{p}}_{p,1}}+\|\dot{S}_{k_0+N_0}f\|_{L^{p^*}}\big)\|g^h\|_{\dot{B}^{\f{n}{p}-1}_{p,1}}\nn\\
&\ls& (\|f^\ell\|_{\dot{B}^{\f{n}{p}}_{2,1}}+\|f\|_{\dot{B}^{1-\f{n}{p}}_{p,1}})\|g^h\|_{\dot{B}^{\f{n}{p}-1}_{p,1}},
\eea where $\f{1}{p*}=\f{1}{2}-\f{1}{p}$. Taking $f=v$ and $g=\nabla O$, and then using \eqref{3.12}, \eqref{3.14} as well as the definition of $\mathcal{G}_p(t)$, we arrive at
\bea\label{3.29}
&&\int^t_0\langle t-\tau\rangle^{-\f{s+s_0}{2}}\|(v\cdot\nabla O^h)\|^\ell_{\dot{B}^{-s_0}_{2,\i}}d\tau\nn\\
&\ls&\int^t_0\langle t-\tau\rangle^{-\f{s+s_0}{2}}(\|v^\ell\|_{\dot{B}^{\f{n}{p}}_{2,1}}+\|v\|_{\dot{B}^{1-\f{n}{p}}_{p,1}})\|\nabla O^h\|_{\dot{B}^{\f{n}{p}-1}_{p,1}}d\tau\nn\\
&\ls&\mathcal{G}^2_p(t)\int^t_0\langle t-\tau\rangle^{-\f{s+s_0}{2}}\big(\langle\tau\rangle^{-(\f{3n}{2p}-\f{n}{4})}+\langle\tau\rangle^{-\f{1}{2}}+\tau^{-\f{1}{2}}
\langle\tau\rangle^{-\f{\al}{2}}\big)\langle \tau\rangle^{-\al}d\tau\nn\\
&\ls&\langle t\rangle^{-\f{s+s_0}{2}}\mathcal{G}^2_p(t).
\eea
Next by taking $f=O$ and $g=\nabla a$ in \eqref{3.28},  we get
\bea\label{3.30}
&&\int^t_0\langle t-\tau\rangle^{-\f{s+s_0}{2}}\|O\cdot\na a^h\|^\ell_{\dot{B}^{-s_0}_{2,\i}}d\tau\nn\\
&\ls&\int^t_0\langle t-\tau\rangle^{-\f{s+s_0}{2}}(\|O^\ell\|_{\dot{B}^{\f{n}{p}}_{2,1}}+\|O\|_{\dot{B}^{1-\f{n}{p}}_{p,1}})\|\nabla a^h\|_{\dot{B}^{\f{n}{p}-1}_{p,1}}d\tau.
\eea
It follows from \eqref{3.12} and \eqref{3.13} that
\bea\label{3.31}
\|O^\ell\|_{\dot{B}^{\f{n}{p}}_{2,1}}+\|O\|_{\dot{B}^{1-\f{n}{p}}_{p,1}}
&\ls&(\langle\tau\rangle^{-(\f{3n}{2p}-\f{n}{4})}+\langle\tau\rangle^{-\f{1}{2}}+\langle\tau\rangle^{-\alpha})\mathcal{G}_p(\tau).
\eea
Consequently, we can get
\bea\label{3.32}
&&\int^t_0\langle t-\tau\rangle^{-\f{s+s_0}{2}}\|O\na a^h\|^\ell_{\dot{B}^{-s_0}_{2,\i}}d\tau\nn\\
&\ls&\mathcal{G}^2_p(t)\int^t_0\langle t-\tau\rangle^{-\f{s+s_0}{2}}(\langle\tau\rangle^{-(\f{3n}{2p}-\f{n}{4})}+\langle\tau\rangle^{-\f{1}{2}}+\langle\tau\rangle^{-\alpha})\langle \tau\rangle^{-\al}d\tau\nn\\
&\ls&\mathcal{G}^2_p(t)\langle t\rangle^{-\f{s+s_0}{2}}.
\eea
To bound the last term $O\nabla v^h$, we need to take $f=O$ and $g=\nabla v$ in \eqref{3.28} and obtain
\bea\label{3.33}
&&\int^t_0\langle t-\tau\rangle^{-\f{s+s_0}{2}}\|O\nabla v^h\|^\ell_{\dot{B}^{-s_0}_{2,\i}}d\tau\nn\\
&\ls&\int^t_0\langle t-\tau\rangle^{-\f{s+s_0}{2}}(\|O^\ell\|_{\dot{B}^{\f{n}{p}}_{2,1}}+\|O\|_{\dot{B}^{1-\f{n}{p}}_{p,1}})\|\nabla v^h\|_{\dot{B}^{\f{n}{p}-1}_{p,1}}d\tau.
\eea
By interpolation, for all $\tau\geq0$,
\be
\|\nabla v^h\|_{\dot{B}^{\f{n}{p}-1}_{p,1}}\ls \Big(\| v\|^h_{\dot{B}^{\f{n}{p}-1}_{p,1}}\|\nabla v\|^h_{\dot{B}^{\f{n}{p}}_{p,1}}\Big)^{\f{1}{2}}\ls\tau^{-\f{1}{2}}\langle\tau\rangle^{-\f{\al}{2}}\mathcal{G}_p(\tau).\nn
\ee
Therefore, we are led to
\bea\label{3.34}
&&\int^t_0\langle t-\tau\rangle^{-\f{s+s_0}{2}}\|O\nabla v^h\|^\ell_{\dot{B}^{-s_0}_{2,\i}}d\tau\nn\\
&\ls&\mathcal{G}^2_p(t)\int^t_0\langle t-\tau\rangle^{-\f{s+s_0}{2}}\langle \tau\rangle^{-\min(\f{1}{2},\f{3n}{2p}-\f{n}{4},\alpha)}\tau^{-\f{1}{2}}\langle\tau\rangle^{-\f{\al}{2}}d\tau\nn\\
&\ls&\mathcal{G}^2_p(t)\langle t\rangle^{-\f{s+s_0}{2}}.
\eea

Putting together all the above estimates for those terms involving in $O$ and those computations with respect to $a$ and $v$  (see \cite{DX1}), we finish the proof of \eqref{3.9}. Then by combining \eqref{6.3} and \eqref{3.9}, we deduce that
\be\label{3.35}
\langle t\rangle^{-\f{s+s_0}{2}}\|(a,O;v)(t)\|^\ell_{\dB^{s}_{2,1}}\ls \mathcal{G}_{p,0}+\mathcal{G}^2_p(t)+\mathcal{X}^2_p(t).
\ee

\subsection{Decay estimates for the high frequencies of $(\nabla a,\nabla O; v)$}
This part is devoted to estimating the second term in $\mathcal{G}_p(t)$. The usual Duhamel principle cannot remain true, since there is a loss of one derivative for the density and deformation tensor at high frequencies. To eliminate the technical difficulty, as the proof of global-in-time existence, we need to perform a suitable quasi-diagonalization, say \textit{effective velocities}, to handle high frequencies.
Let $\dD_k w=w_k$ and $\dD_k \O^{ij}=\O^{ij}_k$. From \eqref{e4.19} and \eqref{ee4.24}, by employing the energy methods of $L^p$ type, we obtain
\bea\label{ee6.35}
\f{d}{dt}\|w_k\|^p_p+c_p2^{2k}\|w_k\|^p_p\ls \big\{2^{-k}\|a_k\|_p+2^{-k}\|\dD_k\t{G}_1\|_p+\|\dD_k G_2\|_p\}\|w_k\|^{p-1}_p,
\eea
\bea\label{ee6.36}
\f{d}{dt}\|\O^{ij}_k\|^p_p+c_p2^{2k}\|\O^{ij}_k\|^p_p\ls \big\{2^{-k}\|O^{ij}_k\|_p+2^{-k}\|\dD_k\t{G}^{ij}_3\|_p+\|\dD_k \t{G}^{ij}_4\|_p\big\}\|\O^{ij}_k\|^{p-1}_p,
\eea
\bea\label{ee6.37}
&&\f{d}{dt}\|\La a_k\|^p_p+c_p\|\La a_k\|^p_p\nonumber\\&\ls& \big\{\|\La\dD_k \t{G}_1\|_p+2^{2k}\|w\|_p+\|\text{div}v\|_\i\|\La a_k\|_p+\|\mathcal{R}_{1,k}\|_p\big\}\|\La a_k\|^{p-1}_p
\eea
and
\bea\label{ee6.38}
&&\f{d}{dt}\|\La O^{ij}_k\|^p_p+c_p\|\La O^{ij}_k\|^p_p\nonumber\\&\ls& \big\{\|\La\dD_k \t{G}^{ij}_3\|_p+2^{2k}\|\O^{ij}\|_p+\|\text{div}v\|_\i\|\La O^{ij}_k\|_p+\|\mathcal{R}_{2,k}\|_p\big\}\|\La O^{ij}_k\|^{p-1}_p
\eea
with $\mathcal{R}_{1,k}\triangleq[v\cdot\na,\La \dD_k]a$ and $\mathcal{R}_{2,k}\triangleq[v\cdot\na,\La \dD_k]O^{ij}$, where we chosen $R_{0}$ sufficiently large such that $2^k>R_0$. Furthermore, with aid of Lemma \ref{lxx}, there exists a constant $c_p>0$ such that
\bea\label{ee6.39}
&&\f{d}{dt}\big(\|\La a_k\|_{\dl,L^p}+\|d_k\|_{\dl,L^p}\big)+c_p\big(\|\La a_k\|_{\dl,L^p}+2^{2k}\|d_k\|_{\dl,L^p}\big)\nn\\
&\ls& \|\La\dD_k \t{G}_1\|_p+\|\dD_k G_2\|_p+\|\text{div}v\|_\i\|\La a_k\|_p+\|\mathcal{R}_{1,k}\|_p+2^{2k}\dl,
\eea
where we used the effective velocity in terms of $a$ and $d$. A similar estimate for $O^{ij}$ and $ e^{ij}$ stems from  \eqref{ee6.36} and \eqref{ee6.38}:
\bea\label{ee6.40}
&&\f{d}{dt}\big(\|\La O_k\|_{\dl,L^p}+\|e^{ij}_k\|_{\dl,L^p}\big)+\tilde{c}_p\big(\|\La O^{ij}_k\|_{\dl,L^p}+2^{2k}\|e^{ij}_k\|_{\dl,L^p}\big)\nn\\
&\ls& \|\La\dD_k \t{G}^{ij}_3\|_p+\|\dD_k \t{G}^{ij}_4\|_p+\|\text{div}v\|_\i\|\La O^{ij}_k\|_p+\|\mathcal{R}_{2,k}\|_p+2^{2k}\dl
\eea
for some constant $\tilde{c}_p>0$. It's easy to see that
\bea\label{ee6.41}
&&\|\dD_k \t{G}^{ij}_4\|_p\ls \|\dD_k {G}^{ij}_4\|_p+2^{2k}\|d_k\|_p+\|\La a_k\|_p\nn\\
&&\qq\qq\ \ls \|\dD_k {G}^{ij}_4\|_p+2^{2k}\|d_k\|_{\dl,L^p}+\|\La a_k\|_{\dl,L^p}.
\eea
Therefore, it follows from \eqref{ee6.39}, \eqref{ee6.40} and \eqref{ee6.41} that
\bea\label{ee6.42}
&&\f{d}{dt}\big(\|(\La a_k,\La O_k; v_k)\|_{\dl,L^p}+c_0\|(\La a_k,\La O_k;v_k)\|_{\dl,L^p}\nn\\
&\ls& \|\text{div}v\|_\i\|(\La a_k,\La O_k)\|_p+\|\mathcal{R}_{1,k}\|_p+\|\mathcal{R}_{2,k}\|_p\nn\\
&&+\|\La\dD_k (\t{G}_1,\t{G}_3)\|_p+\|\dD_k (G_2,G_4)\|_p+2^{2k}\dl
\eea
for some constant $c_0>0$. Integrating in time on both sides and letting $\dl\rightarrow 0$, we eventually get
\be\label{3.44}
\|(\La a_k,\La O_k; v_k)(t) \|_p\leq e^{-c_0t}\|(\La a_k(0),\La O_k(0),v_k(0))\|_p+C\int^t_0 e^{c_0(\tau-t)}g_k(\tau)d\tau,
\ee
where
\bea
g_k&\triangleq& \underbrace{\|\text{div}v\|_{L^\i}\big(\|\La a_k\|_p+\|\La O_k\|_p\big)}_{g^1_k}\nn\\
&&+\underbrace{\|\dD_k ({G}_2, {G}_4,\La \t{G}_1,\La \t{G}_3)\|_p}_{g^2_k}\nn\\
&&+\underbrace{\|\big(\mathcal{R}_{1,k},\mathcal{R}_{2,k}\big)\|_p}_{g^3_k}.\nn
\eea
Multiplying \eqref{3.44} by $\langle t \rangle^{\al}2^{k(\f{n}{p}-1)}$, taking the supremum on $[0,T]$ and summing up over $k$ satisfying $2^k> R_0$ yields
\bea\label{3.45}
\|\langle t\rangle^{\al}(\La a,\La O;v)\|^h_{\t{L}^\i_T\dot{B}^{\f{n}{p}-1}_{p,1}}&\ls&\|(\La a_0,\La O_0;v_0)\|^h_{\dot{B}^{\f{n}{p}-1}_{p,1}}\nonumber\\
&&+\sum\limits_{2^k> R_0} \sup\limits_{0\leq t\leq T}\Bigg(\ltr\int^t_0 e^{c_0(\tau-t)}2^{k(\f{n}{p}-1)}g_k(\tau)d\tau\Bigg).
\eea
 Without loss of generality, we assume that $T\geq 2$ and first bound the the supremum for $0\leq t\leq 2$.
Notice that
\bea\label{3.46}
\sum\limits_{2^k> R_0} \sup\limits_{0\leq t\leq 2}\Bigg(\ltr\int^t_0 e^{c_0(\tau-t)}2^{k(\f{n}{p}-1)}g_k(\tau)d\tau\Bigg)
\ls \int^2_0 \sum\limits_{2^k> R_0}2^{k(\f{n}{p}-1)}g_k(\tau)d\tau.
\eea
Furthermore, it will be shown that the right side of \eqref{3.46} can be bounded by $\mX^2_p(2)$. Indeed,  using Lemma\ref{ll2.6} and the representation of $\t{G}_i(i=1,3)$ and ${G}_i(i=2,4)$, we get
\bea\label{3.47}
&&\int^2_0 \sum\limits_{2^k> R_0}2^{k(\f{n}{p}-1)}g_k(\tau)d\tau\nn\\
&\leq& \int^2_0 \bigg\{\underbrace{\|\text{div} v\|_{L^\i}\big(\|\La O\|^h_{\dot{B}^{\f{n}{p}-1}_{p,1}}+\|\La a\|^h_{\dot{B}^{\f{n}{p}-1}_{p,1}}\big)}_{\text{coming\ from}\ g^1_k}\nn\\
&&+\underbrace{\|(\text{div}(aO),\La(\p_k v^iO^{kj}),O^{lk}\p_lO^{ik},O^{lj}\p_lO^{ik},O^{lk}\p_lO^{ij},O^{jk}\p_jO^{\bullet k})\|^h_{\dot{B}^{\f{n}{p}-1}_{p,1}}}_{\text{coming\ from}\ g^2_k}\nn\\
&&+\underbrace{\Big\|\Big(a\cdot\na v,v\cdot\na v,I(a)\mathcal{A}v,K(a)\na a,\f{1}{1+a}\text{div}\big(2\t{\mu}(a)D(v)+\t{\la}(a)\text{div} v\text{Id}\big)\Big)\Big\|^h_{\dot{B}^{\f{n}{p}-1}_{p,1}}}_{\text{coming\ from}\ g^2_k}\nn\\
&&+\underbrace{\|\na v\|_{\dB^{\f{n}{p}}_{p,1}}\big(\|O\|_{\dB^{\f{n}{p}}_{p,1}}+(\|a\|_{\dB^{\f{n}{p}}_{p,1}}\big)}_{\text{coming\ from}\ g^3_k}\bigg\}d\tau.
\eea
In comparison with \cite{DX1}, we focus on those terms only involving $O$. For instance, we see that
\bea
\int^2_0 \Big(\|\text{div} v\|_{L^\i}\|\La O\|^h_{\dot{B}^{\f{n}{p}-1}_{p,1}}+\|\na v\|_{\dB^{\f{n}{p}}_{p,1}}\|O\|_{\dB^{\f{n}{p}}_{p,1}}\Big)d\tau
\ls \|O\|_{L^\i_t\dot{B}^{\f{n}{p}}_{p,1}}\int^2_0 \| v\|_{\dot{B}^{\f{n}{p}+1}_{p,1}}d\tau
\ls\mX^{2}_p(2).\nn
\eea
Owing to Lemma \ref{ll2.3} and  interpolation  inequalities, we obtain
\bea
&&\int^2_0 \|(\text{div}(aO),\La(\p_k v^iO^{kj}),O^{lk}\p_lO^{ik},O^{lj}\p_lO^{ik},O^{lk}\p_lO^{ij},O^{jk}\p_jO^{\bullet k}\|^h_{\dot{B}^{\f{n}{p}-1}_{p,1}}d\tau\nn\\
&\ls& \int^2_0 \|a\|_{\dot{B}^{\f{n}{p}}_{p,1}}\|O\|_{\dot{B}^{\f{n}{p}}_{p,1}}+\|v\|_{\dot{B}^{\f{n}{p}+1}_{p,1}}\|O\|_{\dot{B}^{\f{n}{p}}_{p,1}}+\|O\|_{\dot{B}^{\f{n}{p}}_{p,1}}\|\na O\|_{\dot{B}^{\f{n}{p}-1}_{p,1}}d\tau\nn\\
&\ls&\Big(\|a\|_{L^2\dot{B}^{\f{n}{p}}_{p,1}}+\|O\|_{L^2\dot{B}^{\f{n}{p}}_{p,1}}\Big)\|O\|_{L^2\dot{B}^{\f{n}{p}}_{p,1}}+\|v\|_{L^1\dot{B}^{\f{n}{p}+1}_{p,1}}\|O\|_{L^\i\dot{B}^{\f{n}{p}}_{p,1}}.\nn\\
&\ls&\mX^{2}_p(2).\nn
\eea
Hence, we can infer that
\bea\label{3.48}
&&\sum\limits_{2^k>R_0}\sup_{0\leq t\leq 2}\ltr \int^t_0 e^{c_0(\tau-t)}2^{k(\f{n}{p}-1)}g_k(\tau)d\tau\ls \mX^{2}_p(2).
\eea
Let us now bound the supremum for $2\leq t\leq T$ in the last term of \eqref{3.45}. To this end, we split the integral on $[0,t]$ into integrals on $[0,1]$ and $[1,t]$. The integral on $[0,1]$ is easy to handle: because $e^{c_0(\tau-t)}\leq e^{-\f{c_0}{2}t}$ for $2\leq t\leq T$ and $0\leq\tau\leq 1$, so one can write
\bea\label{3.49}
\sum\limits_{2^k> R_0}\sup\limits_{2\leq t\leq T}\ltr \int^1_0 e^{c_0(\tau-t)}2^{k(\f{n}{p}-1)}g_k(\tau)d\tau
&\ls& \sum\limits_{2^k> R_0}\sup\limits_{2\leq t\leq T}\ltr e^{-\f{c_0}{2}t} \int^1_0 2^{k(\f{n}{p}-1)}g_k(\tau)d\tau\nn\\
&\ls&  \sum\limits_{2^k> R_0} \int^1_0 2^{k(\f{n}{p}-1)}g_k(\tau)d\tau.
\eea
Therefore, following the procedure leading to \eqref{3.48}, we end up with
\bea
\sum\limits_{2^k> R_0}\sup\limits_{2\leq t\leq T}\ltr \int^1_0 e^{c_0(\tau-t)}2^{k(\f{n}{p}-1)}g_k(\tau)d\tau \ls \mX^{2}_p(1).
\eea
In order to bound the integral on $[1,t]$ for $2\leq t\leq T$, we notice that
\be\label{3.50}
\sum\limits_{2^k> R_0}\sup\limits_{2\leq t\leq T}\Bigg(\ltr \int^t_1 e^{c_0(\tau-t)}2^{k(\f{n}{p}-1)}g_k(\tau)d\tau\Bigg)\ls \sum\limits_{2^k> R_0}2^{k(\f{n}{p}-1)}\sup\limits_{1\leq t\leq T}(t^{\al} g_k(t)).
\ee
In nonlinear sources $g^1_k,g^2_k$ and $g^3_k$, the calculations for those terms with respect to $O$ are totally similar, so we present them for $\mathrm{div}(aO)$ and $\La(\p_k v^iO^{kj})$ for brevity.
We write \be
\text{div}(aO)=a\na\cdot O+\na a\cdot O. \nn
\ee
Due to the same regularity level, it suffices to estimate the term $a\na\cdot O$.
Using Lemma \ref{ll2.3}, we deduce that
\bea\label{3.56}
\|t^{\al}(a\na\cdot O^h)\|_{\t{L}^\i_T\dot{B}^{\f{n}{p}-1}_{p,1}}
\ls \|a\|_{\t{L}^\i_T\dot{B}^{\f{n}{p}}_{p,1}}\|t^{\al}\na O\|^h_{\t{L}^\i_T\dot{B}^{\f{n}{p}-1}_{p,1}}\leq \mX_p(T)\mathcal{G}_p(T),
\eea
and
\bea\label{3.57}
\hspace{-5mm}\|t^{\al}(a\na\cdot O^\ell)\|_{\t{L}^\i_T\dot{B}^{\f{n}{p}-1}_{p,1}}
&\ls& \|t^{\f{\al}{2}}a\|_{\t{L}^\i_T\dot{B}^{\f{n}{p}}_{p,1}}\|t^{\f{\al}{2}}O\|^\ell_{\t{L}^\i_T\dot{B}^{\f{n}{p}}_{p,1}}\nn\\
&\ls& \Big(\|t^{\f{\al}{2}}a\|^\ell_{\t{L}^\i_T\dot{B}^{\f{n}{2}}_{2,1}}+\|t^{\f{\al}{2}}a\|^h_{\t{L}^\i_T\dot{B}^{\f{n}{p}}_{p,1}}\Big)\|t^{\f{\al}{2}}O\|^\ell_{\t{L}^\i_T\dot{B}^{\f{n}{2}}_{2,1}}\leq \mathcal{G}^2_p(T),
\eea
since the fact $\f{\al}{2}\leq \frac{s_0}{2}+\frac{n}{4}-\frac{\varepsilon}{2} $ indicates that $\|t^{\f{\al}{2}}z\|^\ell_{\t{L}^\i_T\dot{B}^{\f{n}{2}}_{2,1}}\lesssim \|t^{\f{\al}{2}}z\|^\ell_{L^\i_T\dot{B}^{\f{n}{2}-\varepsilon}_{2,1}}\lesssim \mathcal{G}_p(T)$ for $z=a,O,v$. Combining \eqref{3.56} and \eqref{3.57}, we get
\be\label{3.58}
\|t^{\al}(a\na\cdot O)\|_{\t{L}^\i_T\dot{B}^{\f{n}{p}-1}_{p,1}}\ls \mX_p(T)\mathcal{G}_p(T)+\mathcal{G}^2_p(T).
\ee
In addition, it follows from \eqref{dfD} that
\be\label{3.51}
\|\tau \na v\|_{\t{L}^\i_t\dot{B}^{\f{n}{p}}_{p,1}}\ls \mathcal{G}_p(t).
\ee
By Lemmas \ref{lem2.1} and \ref{ll2.3}, we have
\bea
\|t^{\al}\La(\p_k v^iO^{kj})\|^h_{\t{L}^\i_T\dot{B}^{\f{n}{p}-1}_{p,1}}&\ls&\|t^{\al}(O^{kj}\p_k v^i)\|^h_{\t{L}^\i_T\dot{B}^{\f{n}{p}}_{p,1}}\nn\\
&\ls&\|t^{\al-1}O \|_{\t{L}^\i_T\dot{B}^{\f{n}{p}}_{p,1}}\|t\na v\|_{\t{L}^\i_T\dot{B}^{\f{n}{p}}_{p,1}}\nn\\
&\ls&\Big(\|t^{\al-1}O \|^\ell_{\t{L}^\i_T\dot{B}^{\f{n}{p}}_{p,1}}+\|t^{\al-1}O \|^h_{\t{L}^\i_T\dot{B}^{\f{n}{p}}_{p,1}}\Big)\|t\na v\|_{\t{L}^\i_T\dot{B}^{\f{n}{p}}_{p,1}}.
\eea
It is obvious that $\|t^{\al-1}O \|^h_{\t{L}^\i_T\dot{B}^{\f{n}{p}}_{p,1}}\leq \mathcal{G}_p(T)$ according to the definition of $\mathcal{G}_p(T)$. On the other hand, we have the following estimates for $z=a,O,v$,
\be\label{3.54}
\|t^{\al-1}z\|^\ell_{\t{L}^\i_T\dot{B}^{\f{n}{p}-1}_{p,1}}\ls\|t^{\al-1}z \|^\ell_{L^\i_T\dot{B}^{\f{n}{2}-1-2\varepsilon}_{2,1}}\leq \mathcal{G}_p(T),
\ee
as $\al-1=\f{1}{2}(s_0+n/2-1-2\varepsilon)$ with enough small $\varepsilon$. Consequently, we arrive at
\be\label{3.55}
\|t^{\al}\La(\p_k v^iO^{kj})\|^h_{\t{L}^\i_T\dot{B}^{\f{n}{p}-1}_{p,1}}\ls \mathcal{G}^2_p(T).
\ee
In a conclusion, by combining those estimates involving $a$ and $v$ in \cite{DX1}, we can conclude that
\bea\label{3.64}
\|\langle t\rangle^{\al} (\La a,\La O;v)\|^h_{\t{L}^\i_T\dot{B}^{\f{n}{p}-1}_{p,1}}
\ls \| (\La a_0,\La O_0;v_0)\|^h_{\dot{B}^{\f{n}{p}-1}_{p,1}}+\mathcal{G}^2_p(T)+\mX^2_p(T).
\eea

\subsection{Decay and gain of regularity for the high frequencies of $v$}
In order to bound the last term in $\mathcal{G}_p(t)$, it is convenient to rewrite the velocity equation in the following way. First, it follows from
(\ref{r1}) that
\bea
\p_tv-\mathcal{A}v&=&F\nn\\
 &\triangleq&-(1+K(a))\na a-v\cdot\nabla v+\na\cdot O+ O^{jk}\p_jO^{\bullet k}-I(a)\mathcal{A}v\\
 &&+\f{1}{1+a}\text{div}\big(2\t{\mu}(a)D(v)+\t{\la}(a)\text{div} v\text{Id}\big)\nn
\eea
Hence, we have
\be\label{3.65}
\p_t(t\mathcal{A}v)-\mathcal{A}(t\mathcal{A}v)=\mathcal{A}v+t\mathcal{A}F.
\ee
We thus deduce from Proposition \ref{pA.4} and the remark that follows, that
\bea\label{3.66}
\|\tau\na^2 v\|^h_{\t{L}^\i_t\dot{B}^{\f{n}{p}-1}_{p,1}}&\ls& \|\mathcal{A} v\|^h_{\t{L}^1_t\dot{B}^{\f{n}{p}-1}_{p,1}}+\|\tau \mathcal{A}F\|^h_{\t{L}^\i_t\dot{B}^{\f{n}{p}-3}_{p,1}}\nn\\
&\ls& \| v\|^h_{\t{L}^1_t\dot{B}^{\f{n}{p}+1}_{p,1}}+\|\tau F\|^h_{\t{L}^\i_t\dot{B}^{\f{n}{p}-1}_{p,1}}\nn\\
&\ls& \mX_p(t)+\|\tau F\|^h_{\t{L}^\i_t\dot{B}^{\f{n}{p}-1}_{p,1}},
\eea
where we have used the bounds given by Theorem \ref{thglobal2}. Secondly, we turn to bound the norm $\|\tau F\|^h_{\t{L}^\i_t\dot{B}^{\f{n}{p}-1}_{p,1}}$.
Because $\al\geq 1$, we have
\be\label{3.67}
\|\tau (\na a,\na\cdot O)\|^h_{\t{L}^\i_t\dot{B}^{\f{n}{p}-1}_{p,1}}\ls \|\langle\tau\rangle^\al(a,O)\|^h_{\t{L}^\i_t\dot{B}^{\f{n}{p}}_{p,1}}.
\ee
Product and composition estimates indicate that
\bea\label{3.68}
\|\tau(K(a)\na a,O^{jk}\p_jO^{\bullet k})\|_{\t{L}^\i_t\dot{B}^{\f{n}{p}-1}_{p,1}}
\ls \|\tau^{\f{1}{2}}(a,O)\|^2_{\t{L}^\i_t\dot{B}^{\f{n}{p}}_{p,1}}\ls \mathcal{G}^2_p(t).
\eea
Together with those estimates for other nonlinear terms (see \cite{DX1}), we can conclude that
\be\label{6.70}
\|\tau\na v\|^h_{\t{L}^\i_t\dot{B}^{\f{n}{p}}_{p,1}}\ls\mX^2_p(t)+\mathcal{G}^2_p(t)+\|\langle \tau\rangle^\al (a,O)\|^h_{\t{L}^\i_t\dot{B}^{\f{n}{p}}_{p,1}}.
\ee
Finally, bounding the last term on the right-side of \eqref{6.70} according to \eqref{3.64}, and
adding up the obtained inequality to \eqref{3.35} and \eqref{6.70} yields for all $t\geq0$
\bea\label{e6.70}
\mathcal{G}_p(t)&\ls& \mathcal{G}_{p,0}+\|(a_0,O_0;v_0)\|^\ell_{\dot{B}^{\f{n}{2}-1}_{2,1}}+\|(\na a_0,\na O_0;v_0)\|^h_{\dot{B}^{\f{n}{p}-1}_{p,1}}+\mathcal{G}^2_p(t)+\mX^2_p(t)\nn\\
&\ls&  \mathcal{G}_{p,0}+\E^{n/p}_0+\mathcal{G}^2_p(t)+\mX^2_p(t).
\eea
It follows from  Theorem\ref{thglobal2} that $\mathcal{X}_p(t)\leq M\E^{n/p}_0\leq M\eta\ll 1$. On the other hand, $\|(a_0,O_0;v_0)\|^\ell_{\dot{B}^{\f{n}{2}-1}_{2,1}}\lesssim \|(a_0,O_0;v_0)\|^\ell_{\dot{B}^{-s_0}_{2,\infty}}$, so one can conclude that
\eqref{de1} is fulfilled for all time if $\mathcal{G}_{p,0}$ and $\|(\na a_0,\na O_0;v_0)\|^h_{\dot{B}^{\f{n}{p}-1}_{p,1}}$ are small enough.
This finishes the proof of  Theorem \ref{thdecay} eventually.

\section*{Appendix: Some Estimates in the Hybrid Besov Space}
\renewcommand{\theproposition}{A.\arabic{proposition}}
\renewcommand\theequation{A.\arabic{equation}}
\renewcommand\theremark{A.\arabic{remark}}
\begin{proposition}
[\cite{CMZ1}\label{A.1}] Let $s_1,s_2,t_1,t_2,\sigma,\tau\in\bR,\ 2\leq p \leq 4$ and $1\leq r,r_1,r_2,r_3,r_4\leq \i$ with $\f{1}{r}=\f{1}{r_1}+\f{1}{r_2}=\f{1}{r_3}+\f{1}{r_4}$. Then we have the following:

\begin{itemize}
\item If $\sigma,\tau\leq n/p$ and $\sigma+\tau>0$, then
\bea\label{A.1}
&&\sum\limits_{2^k> R_0}2^{k(\sigma+\tau-n/p)}\|\dD_k(fg)\|_{L^r_TL^p}\nn\\
&\leq&C\|f\|_{\t{L}^{r_1}_T\B^{n/2-n/p+\sigma,\sigma}_{2,p}}\|g\|_{\t{L}^{r_2}_T\B^{n/2-n/p+\tau,\tau}_{2,p}}.
\eea

\item If $s_1,s_2\leq n/p$ and $s_1+t_1>n-\f{2n}{p}$ with $s_1+t_1=s_2+t_2$ and $\gamma\in\bR$, then
\bea\label{A.2}
&&\sum\limits_{2^k\leq R_0}2^{k(s_1+t_1-n/2)}\|\dD_k(fg)\|_{L^r_TL^2}\nn\\
&\leq&C\big(\|f\|_{\t{L}^{r_1}_T\B^{s_1,s_1-n/2+n/p}_{2,p}}\|g\|_{\t{L}^{r_2}_T\B^{t_1,t_1-n/2+n/p+\gamma}_{2,p}}\nn\\
&&\qq+\|g\|_{\t{L}^{r_3}_T\B^{s_2,s_2-n/2+n/p}_{2,p}}\|f\|_{\t{L}^{r_4}_T\B^{t_2,t_2-n/2+n/p}_{2,p}}\big).
\eea
\item If $s_1,s_2\leq n/2$ and $s_1+t_1>\f{n}{2}-\f{n}{p}$ with $s_1+t_1=s_2+t_2$, then
\bea\label{A.3}
&&\sum\limits_{k\in\bZ}2^{k(s_1+t_1-n/2)}\|\dD_k(fg)\|_{L^r_TL^2}\nn\\
&\leq&C\big(\|f\|_{\t{L}^{r_1}_T\B^{s_1,s_1-n/2+n/p}_{2,p}}\|g\|_{\t{L}^{r_2}_T\dot{B}^{t_1}_{2,1}}\nn\\
&&\qq+\|g\|_{\t{L}^{r_3}_T\B^{s_2,s_2-n/2+n/p}_{2,p}}\|f\|_{\t{L}^{r_4}_T\dot{B}^{t_2}_{2,1}}\big).
\eea
\end{itemize}
\end{proposition}

\begin{proposition}
[\cite{CMZ1}\label{pa.2}] Let $2\leq p \leq 4$, $s,\sigma>0$, and $s\geq \sigma-n/2+n/p,\ r\geq 1$. Assume that $F\in W^{[s]+2,\i}_{loc}\cap W^{[\sigma]+2,\i}_{loc}$ with $F(0)=0$. Then there haves
\be\label{A.4}
\|F(f)\|_{\t{L}^r_T\B^{s,\sigma}_{2,p}}\leq C(1+\|f\|_{\t{L}^\i_T\B^{n/p,n/p}_{2,p}})^{\max([s],[\sigma])+1}\|f\|_{\t{L}^r_T\B^{s,\sigma}_{2,p}}.
\ee
For any $s>0$ and $p\geq 1$, there haves
\be\label{A.5}
\|F(f)\|_{\t{L}^r_T\dB^{s}_{p,1}}\leq C(1+\|f\|_{{L}^\i_TL^\i})^{[s]+1}\|f\|_{\t{L}^r_T\dB^{s}_{p,1}}.
\ee

\end{proposition}
\begin{proposition}[\cite{Dr2}]\label{pa.3}
Let $s\in(-n\min(1/p,1/p'),1+n/p)$ and $1\leq p,q\leq \i$. Let $v$ be a vector field such that $\na v\in L^1_T\dB^{n/p}_{p,1}$. Assume that $f_0\in\dB^s_{p,q},g\in L^1_T\dB^s_{p,q}$, and $f$ is a solution of the transport equation
\be
\p_tf+v\cdot\na f=g,\q f|_{t=0}=f_0.\nn
\ee
Then for $t\in[0,T]$, there holds
\be
\|f\|_{\t{L}_t\dB^s_{p,q}}\leq \exp\big(C\int^t_0\|\na v(\tau)\|_{\dB^{n/p}_{p,1}}d\tau\big)\big(\|f_0\|_{\dB^s_{p,q}}+\int^t_0\|g(\tau)\|_{\dB^s_{p,q}}d\tau\big).\nn
\ee
\end{proposition}
For the heat equation, one has the following optimal regularity estimate.
\begin{proposition}\label{pA.4}
  Let $p,r\in [1,\i]$, $s\in\bR$, and $1\leq \rho_2\leq \rho_1\leq\i $ Assume that $u_0\in\dB^{s-1}_{p,r}$, $f\in \t{L}^{\rho_2}_T\dB^{s-3+\f{2}{\rho_2}}_{p,r}$. Let $u$ be a solution of the equation
 \be
 \p_tu -\mu \Dl u=f,\q u|_{t=0}=u_0.\nn
 \ee
 Then for $t\in[0,T]$, there holds
 \be\label{A.6}
 \mu^{\f{1}{\rho_1}}\|u\|_{\t{L}^{\rho_1}_T\dB^{s-1+2/\rho_1}_{p,r}}\leq C\big(\|u_0\|_{\dB^{s-1}_{p,r}}+\mu^{1/\rho_2-1}\|f\|_{\t{L}^{\rho_2}_T\dB^{s-3+\f{2}{\rho_2}}_{p,r}}\big).
 \ee
\end{proposition}
\begin{remark}
The estimate \eqref{A.6} is still hold for the following equation
 \be\label{A.7}
 \p_tu -\mu \Dl u-(\la+\mu)\na\text{div}u=f,\q u|_{t=0}=u_0,
 \ee
where $\la$ and $\mu$ are constants such that $\mu>0$ and $\la+\mu>0$(up to the different dependence on the viscous coefficients). Indeed, both $\mathcal{P} u$ and $\mathcal{P}^\perp u$ satisfy the heat equation. We can apply $\mathcal{P}$ and $\mathcal{P}^\perp $ to \eqref{A.7} to get the heat estimate \eqref{A.6}.
\end{remark}

\noindent
{\bf Acknowledgement.}
 The first author is supported by the Research Foundation of Nanjing University of Aeronautics and Astronautics(NO. 1008-YAH17070). The second author is partially supported by the National Natural
Science Foundation of China (11471158) and the Fundamental Research Funds for the Central Universities
(NE2015005). He would like to thank Professor Rapha\"el Danchin for his kind communication when visiting the LAMA in UPEC.

\qq \qq\qq\qq\qq\qq Xinghong Pan

\qq\qq\qq\qq\qq\qq Department of Mathematics,

\qq\qq\qq\qq\qq\qq Nanjing University of Aeronautics and Astronautics,

\qq\qq\qq\qq\qq\qq Nanjing 211106, People's Republic of China.

\qq\qq\qq\qq\qq\qq Email:xinghong\underline{\ \ }87@nuaa.edu.cn

\vskip 0.5cm

\qq \qq\qq\qq\qq\qq Jiang Xu

\qq\qq\qq\qq\qq\qq Department of Mathematics,

\qq\qq\qq\qq\qq\qq Nanjing University of Aeronautics and Astronautics,

\qq\qq\qq\qq\qq\qq Nanjing 211106, People's Republic of China.

\qq\qq\qq\qq\qq\qq Email:jiangxu\underline{\ \ }79@nuaa.edu.cn, jiangxu\underline{\ \ }79math@yahoo.com

\end{document}